\documentclass[a4,10pt,amssymb,leqno]{amsart}

\usepackage{graphics}
\usepackage{epsfig}
\usepackage{enumerate}
\usepackage[french]{babel}
\usepackage[T1]{fontenc}

\title[
On wedge extendability of CR-meromorphic functions]
{On wedge extendability of CR-meromorphic functions}

\newtheorem{thm}{Theorem}[section]

\newtheorem{lem}[thm]{Lemma}
\newtheorem{mainlem}[thm]{Main Lemma}

\def\N{{\Bbb N}}

\def\Z{{\Bbb Z}}

\def\R{{\Bbb R}}
\def\C{{\Bbb C}}

\def\dim{\hbox{dim\,}}
\def\codim{\hbox{codim}}
\def\v{\vert}

\def\eps{\varepsilon}

\author{Jo\"el Merker and Egmont Porten}

\address{LATP, UMR 6632,
Centre de Math\'ematiques et d'Informatique,
39 rue Joliot-Curie,
F-13453 Marseille Cedex 13, France}

\email{merker\@@cmi.univ-mrs.fr}

\address{Humboldt-Universit\"at zu Berlin, 
Mathematisch-Naturwissenschaftenschaft\-liche Fa\-kult\"at II, 
Institut f\"ur Mathematik, Rudower Chaussee 25, D-12489 Berlin, Germany}

\email{egmont\@@mathematik.hu-berlin.de}

\keywords{CR-meromorphic functions, Removable singularities for
the tangential Cauchy-Riemann operator, Analytic discs.}

\subjclass{Primary: 32D20. Secondary: 32A20, 32D10, 32V10, 32V25, 32V35}
\date{\number\year-\number\month-\number\day}

\begin{document}


\begin{abstract}
In this article, we consider metrically thin singularities $E$ of the
solutions of the tangential Cauchy-Riemann operators on a
$\mathcal{ C}^{2,\alpha}$-smooth 
embedded Cauchy-Riemann generic manifold $M$ (CR
functions on $M\backslash E$) and more generally, we consider
holomorphic functions defined in wedgelike domains attached to
$M\backslash E$. Our main result establishes the wedge- and the
$L^1$-removability of $E$ under the hypothesis that the $(\dim
M-2)$-dimensional Hausdorff volume of $E$ is zero and that $M$ and
$M\backslash E$ are globally minimal. As an application, we deduce
that there exists a wedgelike domain attached to 
an everywhere locally minimal $M$ to which every
CR-meromorphic function on
$M$ extends meromorphically.
\end{abstract}

\maketitle

\section{Introduction and statement of results}

In continuation with our previous works [MP1,2,3], we study the wedge
removability of metrically thin singularities of CR functions and its
application to the local extendability of CR-meromorphic functions 
defined on CR manifolds of arbitrary codimension.

First we need to recall some fundamental notions concerning CR
manifolds. For a detailed presentation we refer to [Bo].  Let $M$ be a
connected smooth CR generic manifold in $\C^{m+n}$ with ${\rm CRdim}\,
M=m\geq 1$, $\codim_{\R}\, M=n\geq 1$, and $\dim_{\R} M=2m+n$. We
denote sometimes $N:=m+n$. In suitable holomorphic coordinates
$(w,z=x+iy)\in\C^{m+n}$, $M$ may be represented as the graph of a
differentiable vector-valued mapping in the form $x=h(w,y)$ with
$h(0)=0$, $d h(0)=0$. The manifold $M$ is called {\it globally
minimal} if it consists of a single CR orbit. This notion generalizes
the concept of {\it local minimality} in the sense of Tumanov, {\it
cf.} [Trp], [Tu1,2], [J1,2], [M], [MP1].  A {\it wedge} $\mathcal{ W}$
{\it with edge} $M'\subset M$ is a set of the form $\mathcal{
W}=\{p+c:p\in M',c\in C\}$, where $C\subset\C^{m+n}$ is a truncated
open cone with vertex in the origin. By a {\it wedgelike domain}
$\mathcal{ W}$ attached to $M$ we mean a domain which contains for
every point $p\in M$ a wedge with edge a neighborhood of $p$ in $M$
({\it cf.} [MP1,2,3]).

A closed subset $E$ of $M$ is called {\it wedge removable} (briefly
$\mathcal{ W}$-removable) if for every wedgelike domain $\mathcal{ W}_1$
attached to $M\backslash E$, there is a wedgelike domain $\mathcal{ W}_2$
attached to $M$ such that for every holomorphic function $f\in\mathcal{
O}(\mathcal{ W}_1)$, there exists a holomorphic function $F\in\mathcal{
O}(\mathcal{ W}_2)$ which coincides with $f$ in some wedgelike open set
$\mathcal{ W}_3\subset \mathcal{ W}_1$ attached to $M\backslash E$.  We say
that $E$ is {\it $L^1$-removable} if every locally integrable function
$f$ on $M$ which is CR on $M\backslash E$ is CR on all of $M$ (here,
CR is understood in the distributional sense).

Let $H^\kappa$ denote $\kappa$-dimensional Hausdorff measure,
$\kappa\geq 0$. Our main result is~:

\begin{thm}\label{main}
Suppose $M$ is $\mathcal{ C}^{2,\alpha}$-smooth, $0<\alpha<1$. Then every
closed subset $E$ of $M$ such that $M$ and $M\backslash E$ are
globally minimal and such that $H^{2m+n-2}(E)=0$ is $\mathcal{ W}$- and
$L^1$-removable.
\end{thm}

(We shall say sometimes that $E$ is of codimension $2^{+0}$ in $M$.)
The hypersurface case of this statement follows from works of
Lupacciolu, Stout, Chirka and others, with weaker regularity
assumptions, $M$ being $\mathcal{ C}^2$-smooth, $\mathcal{ C}^1$-smooth or 
even a Lipschitz graph ({\it see} [LS], [CS]), so Theorem 1.1 is new
essentially in codimension $n\geq 2$. Recently, many geometrical
removability results have been established in case the singularity $E$
is a submanifold ({\it see} [St], [LS], [CS], [J2,3], [P1], [MP1,2,3],
[JS], [P2], [MP4]) and Theorem 1.1 appears to answer one of the last
open general questions in the subject ({\it see} also [J3], [MP4] for
related open problems). As a rule $L^1$-removability follows once
$\mathcal{ W}$-removability being established ({\it see} especially
Proposition 2.11 in [MP1]). In the case at hand we have already proved
$L^1$-removability by different methods earlier (Theorem 3.1 in [MP3])
and also $\mathcal{ W}$-removability if $M$ is real analytic ({\it
see} [MP2, Theorem 5.1], with $M$ being $\mathcal{ C}^\omega$-smooth and
$H^{2m+n-2}(E)=0$).

For the special case where $M$ is $\mathcal{ C}^3$-smooth and
Levi-nondegenerate ({\it i.e.}~the convex hull of the image of the
Levi-form has nonempty interior), Theorem~1.1 is due to Dinh and
Sarkis [DS].  It is known that this assumption entails the dimensional
inequality $m^2\geq n$.  Especially, in the case of CR dimension
$m=1$, the abovementioned authors recover only the known hypersurface
case ($n=1$). We also point out a general restriction: by assuming that
$M$ is Levi-nondegenerate, or more generally that it is of
Bloom-Graham finite type at {\it every} point of $M$, one would not
take account of propagation aspects for the regularity of CR
functions. For instance, it is
well known that wedge extendability may hold despite of
large Levi-flat regions in manifolds $M$ consisting of a single 
CR orbit ({\it cf.} [Trp], [Tu1,2], [J1], [M]). For the sake of
generality, this is why we only assume that $M$ and $M\backslash E$
are globally minimal in Theorem~1.1.

A straightforward application is as follows. First, by [Trp], [Tu1,2],
[J1], [M, Theorem 3.4], as $M\backslash E$ is globally minimal, there
is a wedgelike domain $\mathcal{ W}_0$ attached to $M\backslash E$ to
which every continuous CR function (resp.  CR distribution) $f$ on
$M\backslash E$ extends as a holomorphic function with continuous
(resp.  distributional) boundary value $f$. Then Theorem~1.1 entails
that there exists a wedge $\mathcal{ W}$ attached to $M$ such that
every such $f$ extend holomorphically as an $F\in \mathcal{ O}(\mathcal{
W})$. There is {\it a priori} no growth control of $F$ up to
$E$. However, as proved in [MP1, Proposition~2.11], in the case where $f$ is
an element of $L^1(M)$ which is CR on $M\backslash E$, some growth
control of Hardy-spaces type can be achieved on $F$ to show that it
admits a boundary value $b(F)$ over $M$ (including $E$) which is $L^1$
and CR on $M$. This is how one may deduce $L^1$-removability from $\mathcal{
W}$-removability in Theorem~1.1.

We now indicate a second application of Theorem 1.1 to the extension
of CR-meromorphic functions. This notion was introduced for
hypersurfaces by Harvey and Lawson [HL] and for generic CR manifolds by
Dinh and Sarkis.  Let $f$ be a CR-meromorphic function, namely: {\bf
1.} $f: \mathcal{ D}_f \to P_1(\C)$ is a $\mathcal{ C}^1$-smooth mapping
defined over a dense open subset $\mathcal{ D}_f$ of $M$ with values in
the Riemann sphere; {\bf 2.} The closure $\Gamma_f$ of its graph in
$\C^{m+n} \times P_1(\C)$ defines an oriented scarred $\mathcal{
C}^1$-smooth CR manifold of CR dimension $m$ ({\em i.e.} CR outside a
closed thin set) and {\bf 3.} We assume that $d[\Gamma_f]=0$ in the
sense of currents ({\it see} [HL], [Sa], [DS], [MP2] for further
definition).  According to an observation of Sarkis based on a
counting dimension argument, the indeterminacy set $\Sigma_f$ of $f$
is a closed subset of empty interior in a two-codimensional scarred
submanifold of $M$ and its scar set is always metrically thin~:
$H^{2m+n-2}(Sc(\Sigma_f))=0$. Moreover, outside $\Sigma_f$, $f$
defines a CR current in some suitable projective chart, hence it
enjoys all the extendability properties of an usual CR function or
distribution. However, the complement $M\backslash \Sigma_f$ need not
be globally minimal if $M$ is, and it is easy to construct manifolds
$M$ and closed sets $E\subset M$ with $H^{2m-1}(E)<\infty$
($m=\dim_{CR} M$) which perturb global minimality ({\it see} [MP1],
p.~811). It is therefore natural to make the additional assumption
that $M$ is minimal (locally, in the sense of Tumanov) at {\it every} point,
which seems to be the weakest assumption which insures that
$M\backslash E$ is globally minimal for arbitrary closed sets
$E\subset M$ (even with a bound on their Hausdorff
dimension). Finally, under these circumstances, the set $\Sigma_f$
will be $\mathcal{ W}$-removable: for its regular part
$Reg(\Sigma_f)$, this already follows from Theorem 4 (ii) in [MP1] and
for its scar set $Sc(\Sigma_f)$, this follows from Theorem 1.1
above. The removability of $\Sigma_f$ means that the envelope of
holomorphy of every wedge $\mathcal{ W}_1$ attached to $M\backslash
\Sigma_f$ contains a wedge $\mathcal{ W}_2$ attached to $M$. As
envelopes of meromorphy and envelopes of holomorphy of domains in
$\C^{m+n}$ coincide by a theorem of Ivashkovich ([I]), we conclude~:

\begin{thm}
Suppose $M$ is $\mathcal{ C}^{2,\alpha}$-smooth and locally minimal at
every point. Then there exists a wedgelike domain $\mathcal{ W}$
attached to $M$ to which all CR-mero\-morphic functions on $M$ extend
meromorphically.
\end{thm}

The remainder of the paper is devoted to the proof of Theorem~1.1.  We
combine the {\it local} and the {\it global} techniques of
deformations of analytic discs, using in an essential way two
important papers of Tumanov [Tu1] and of Globevnik [G1]. In Sections~2
and 3, we first set up a standard local situation ({\it cf.}
[MP1,2,3]).  These preliminaries provide the necessary background for
an informal discussion of the techniques of deformations of analytic
discs we have to introduce. After these motivating remarks, a detailed
presentation of the main part of the proof is provided in Section~4
({\it see} especially Main Lemma~4.3).

\smallskip
\noindent
{\sf Acknowledgements.} The authors would like to thank N. Eisen and
F. Sarkis for instructing discussions in this subject. This article
was conceived when the second author was visiting the LATP, UMR 6632
du CNRS, at the University of Provence, Aix-Marseille I. The second
author is very grateful to this institution for its hospitality.

\section{Localization}

The following section contains important preliminary steps for the
proof of Theorem 1.1 ({\it cf.} [MP1,2,3]). 

As in [CS, p.96], we shall proceed by contradiction, since this
strategy simplifies the general reasonings in the large. Also, 
in Section~3 below, we shall explain how to reduce the question to the
simpler case where the functions which we have to extend are even
holomorphic near $M\backslash E$. Whereas such a
strategy is carried out in detail in [MP1] (with minor variations), we
shall for completeness recall the complete reasonings briefly here, in
Sections~2 and~3.

Thus, we fix $\mathcal{ W}_1$ attached to $M\backslash E$ and say that
an open submanifold $M'\subset M$ containing $M\backslash E$ enjoys
the {\it $\mathcal{ W}$-extension property} if there is a wedgelike
domain $\mathcal{ W}'$ attached to $M'$ and a wedgelike set $\mathcal{
W}_1'\subset\mathcal{ W'}\cap\mathcal{ W}_1$ attached to $M\backslash
E$ such that, for each function $f\in\mathcal{ O}(\mathcal{ W}_1)$,
its restriction to $\mathcal{ W}'_1$ extends holomorphically to
$\mathcal{ W}'$. 

This notion can be localized as follows. Let $E'\subset E$
be an arbitrary closed subset of $E$. We shall say that a point
$p'\in E'$ is (locally) {\it removable} (with respect to $E'$) if for
every wedgelike domain $\mathcal{ W}_1$ attached to $M\backslash E'$,
there exists a neighborhood $U$ of $p'$ in $M$ and a wedgelike domain
$\mathcal{ W}_2$ attached to $(M\backslash E') \cup U$ such that for
every holomorphic function $f\in \mathcal{ O}(\mathcal{ W}_1)$, there
exists a holomorphic function $F\in \mathcal{ O}(\mathcal{ W}_2)$
which coincides with $f$ in some wedgelike open set $\mathcal{
W}_3\subset \mathcal{ W}_1$ attached to $M\backslash E'$.

Next, we define the following set of closed subsets of $E$:
$$
\aligned
\mathcal{ E}:=\{E'\subset E \ \hbox{closed} \ ;
\ M\backslash E' \ & \text{\rm is globally minimal}\\
&\text{\rm and has the 
$\mathcal{ W}$-extension property} \}.
\endaligned
$$
Then the residual set 
$$
E_{\text{\rm nr}}:= \bigcap_{E' \in\mathcal{ E}} E'
$$ 
is closed.  Here, the letters ``nr'' abbreviate ``non-remova\-ble'',
since one expects {\it a priori} that no point of $E_{\rm nr}$ should
be removable in the above sense.  Notice that for any two sets
$E_1',E_2'\in\mathcal{ E}$, $M\backslash E_1'$ and $M\backslash E_2'$
consist of a single CR orbit and have nonempty intersection. Hence
$(M\backslash E_1') \cup (M\backslash E_2')$ is globally minimal and
it follows that $M\backslash E_{\text{\rm nr}}$ is globally minimal.

Using Ayrapetian's version of the edge of the wedge theorem ({\it see}
also [Tu1, Theorem~1.2]), the different wedgelike domains attached to
the sets $M\backslash E'$ can be glued (after appropriate contraction
of their cone) to a wedgelike domain $\mathcal{ W}_1$ 
attached to $M\backslash E_{\text{\rm
nr}}$ in such a way that $M\backslash E_{\text{\rm nr}}$ enjoys the
$\mathcal{ W}$-extension property.  Clearly, to establish Theorem~1.1,
it is enough to show that $E_{\text{\rm nr}}=\emptyset$.

Let us argue indirectly (by contradiction) 
and assume that $E_{\text{\rm nr}}\neq
\emptyset$. With respect to the ordering of sets by the inclusion
relation, $E_{\rm nr}$ is then the minimal non-removable subset of $E$. In
order to derive a contradiction to the minimality of $E_{\text{\rm
nr}}$, it suffices therefore to remove one single point $p\in E_{\text{\rm
nr}}$.  More precisely one has to look for a neighborhood $U_{p}$ of
$p$ such that $U_{p}\cup (M\backslash E_{\text{\rm nr}})$ is
globally minimal and has the $\mathcal{ W}$-extension property.

In order to achieve the first required property, it is very convenient
to choose the point $p$ such that locally the singularity $E_{\rm nr}$
lies behind a ``wall'' through $p$. More precisely we shall construct
a generic hypersurface $M_1\subset M$ containing $p$ such that a
neighborhood $V$ of $p$ in $M$ writes as the disjoint union $M^+\cup
M^- \cup M_1$ of connected sets, where $M^\pm$ are two open
``sides'', and the inclusion $E_{\text{\rm nr}}\cap V\subset
M^-\cup\{p\}$ holds true.  Since $M_1$ is a generic CR manifold, there
is a CR vector field $X$ on $M$ defined in a neighborhood of $p$ which
is transverse to $M_1$. By integrating $X$, one easily finds a basis of
neighborhoods $U$ of $p$ in $M$ such that $U\cup (M\backslash
E_{\text{\rm nr}})$ is globally minimal. Hence it remains to establish
the $\mathcal{ W}$-extension property at $p$, which is the main task.

For sake of completeness, we recall from [MP1] how to construct 
the generic wall $M_1$.

\begin{lem}
There is a point $p_1\in E_{\text{\rm nr}}$ and a $\mathcal{
C}^{2,\alpha}$-smooth generic hypersurface $M_1\subset M$ passing
through $p_1$ so that $E_{\text{\rm nr}}\backslash\{p_1\}$ lies near $p_1$
on one side of $M_1$ $(${\rm see}
{\sc Figure~1}$)$.
\end{lem}

{\it Proof.}  Let $p\in E_{\rm nr}\neq \emptyset$ be an arbitrary
point and let $\gamma$ be a piecewise differentiable CR-curve linking
$p$ with a point $q\in M\backslash E_{\rm nr}$ (such a $\gamma$ exists
because $M$ and $M\backslash E_{\rm nr}$ are globally minimal by
assumption). After shortening $\gamma$, we may suppose that
$\{p\}=E_{\rm nr}\cap \gamma$ and that $\gamma$ is a smoothly embedded
segment. Therefore $\gamma$ can be described as a part of an integral
curve of some nonvanishing $\mathcal{ C}^{1,\alpha}$-smooth CR vector field
(section of $T^cM$) $L$ defined in a neighborhood of $p$.

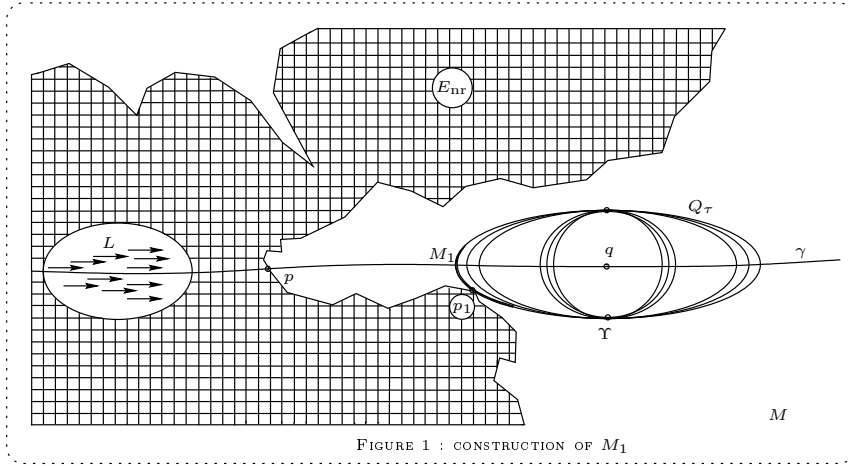
\bigskip
\begin{center}
\begin{picture}(0,0)%
\epsfig{file=figure1.pstex}%
\end{picture}%
\setlength{\unitlength}{3947sp}%
\begingroup\makeatletter\ifx\SetFigFont\undefined%
\gdef\SetFigFont#1#2#3#4#5{%
  \reset@font\fontsize{#1}{#2pt}%
  \fontfamily{#3}\fontseries{#4}\fontshape{#5}%
  \selectfont}%
\fi\endgroup%
\begin{picture}(5356,2919)(1333,-3388)
\put(1959,-2018){\makebox(0,0)[lb]{\smash{\SetFigFont{6}{7.2}{\familydefault}{\mddefault}{\updefault}$L$}}}
\put(6142,-3096){\makebox(0,0)[lb]{\smash{\SetFigFont{6}{7.2}{\familydefault}{\mddefault}{\updefault}$M$}}}
\put(5114,-2054){\makebox(0,0)[lb]{\smash{\SetFigFont{6}{7.2}{\familydefault}{\mddefault}{\updefault}$q$}}}
\put(5070,-2586){\makebox(0,0)[lb]{\smash{\SetFigFont{6}{7.2}{\familydefault}{\mddefault}{\updefault}$\Upsilon$}}}
\put(6312,-2061){\makebox(0,0)[lb]{\smash{\SetFigFont{6}{7.2}{\familydefault}{\mddefault}{\updefault}$\gamma$}}}
\put(3547,-3284){\makebox(0,0)[lb]{\smash{\SetFigFont{6}{7.2}{\familydefault}{\mddefault}{\updefault}{\sc Figure 1: construction of $M_1$}}}}
\put(5639,-1783){\makebox(0,0)[lb]{\smash{\SetFigFont{6}{7.2}{\familydefault}{\mddefault}{\updefault}$Q_\tau$}}}
\put(3100,-2225){\makebox(0,0)[lb]{\smash{\SetFigFont{6}{7.2}{\familydefault}{\mddefault}{\updefault}$p$}}}
\put(4054,-1032){\makebox(0,0)[lb]{\smash{\SetFigFont{6}{7.2}{\familydefault}{\mddefault}{\updefault}$E_{\rm nr}$}}}
\put(4007,-2082){\makebox(0,0)[lb]{\smash{\SetFigFont{6}{7.2}{\familydefault}{\mddefault}{\updefault}$M_1$}}}
\put(4161,-2403){\makebox(0,0)[lb]{\smash{\SetFigFont{6}{7.2}{\familydefault}{\mddefault}{\updefault}$p_1$}}}
\end{picture}

\end{center}
\bigskip

Let $H\subset M$ be a small $(\dim \,M-1)$-dimensional hypersurface of
class $\mathcal{ C}^{2,\alpha}$ passing through $p$ and transverse to
$L$. Integrating $L$ with initial values in $H$ we obtain $\mathcal{
C}^{1,\alpha}$-smooth coordinates $(t,s)\in\R\times\R^{{\rm dim}\, M-1}$ so
that for fixed $s_0$, the segments $(t,s_0)$ are contained in the
trajectories of $L$. After a translation, we may assume that $(0,0)$
corresponds to a point of $\gamma$ close to $p$ which is not contained in
$E_{\text{\rm nr}}$, again denoted by $q$. 
Fix a small $\eps>0$ and for real $\tau\geq 1$, 
define the ellipsoids ({\it see} {\sc Figure~1} above)
$$
Q_\tau:=\{(t,s):|t|^2/\tau+|s|^2<\eps\}.
$$
There is a minimal $\tau_1>1$ with $\overline{Q_{\tau_1}}\cap
E_{\text{\rm nr}}\not=\emptyset$. Then $ \overline{Q_{\tau_1}}\cap
E_{\text{\rm nr}}=\partial Q_{\tau_1}\cap E_{\text{\rm nr}}$ and 
${Q_{\tau_1}}\cap E_{\text{\rm nr}}=\emptyset$.  Observe
that every $\partial Q_{\tau}$ is transverse to the trajectories of $L$ out
off the equatorial set $\Upsilon:=\{(0,s):|s|^2=\eps\}$ which
is contained in 
$M\backslash E_{\text{\rm nr}}$.  Hence $\partial Q_{\tau_1}$ is
transverse to $L$ in all points of $\partial Q_{\tau_1}\cap
E_{\text{\rm nr}}$. So $\partial Q_{\tau_1}\backslash \Upsilon$ is
generic in $\C^{m+n}$, since $L$ is a CR field.

We could for instance choose a point $p_1\in\partial Q_{\tau_1}\cap
E_{\text{\rm nr}}$ and take for $M_1$ a neighborhood of $p_1$ in
$\partial Q_{\tau_1}$, but such an $M_1$ would be {\it only} of class
$\mathcal{ C}^{1,\alpha}$ and we want 
$\mathcal{ C}^{2,\alpha}$-smoothness.

Therefore we fix a small $\delta>0$ and approximate the family
$\partial Q_\tau,1\leq\tau<\tau_1+\delta$, by a nearby family of
$\mathcal{ C}^{2,\alpha}$-smooth hypersurfaces $\partial
\widetilde{Q}_\tau,1\leq \tau<\tau_1+\delta$.  Clearly this can be done so
that the $\partial \widetilde{Q}_\tau$ are still
boundaries of increasing domains $\widetilde{Q}_\tau$ of approximately
the same size as $Q_\tau$ and so that the points where the $\partial
\widetilde{Q}_\tau$ are tangent to $L$ are also
contained in $M\backslash E_{\text{\rm nr}}$ near the equator $\Upsilon$
of $Q_\tau$.

The same reasoning as above shows that there exist a real number
$\widetilde{\tau}_1>1$, a point $p_1\in E_{\rm nr}$ and a generic
hypersurface $M_0$ passing through $p_1$ (which is a piece of
$\partial \widetilde{Q}_{\widetilde{\tau}_1}$) such that $E_{\rm nr}$
lies in the left closed side $M_0^-\cup M_0$ in a neighborhood of
$p_1$ ({\it see} {\sc Figure~1}).  We want more: $E_{\rm nr}\backslash
\{p_1\}\subset M_1^-$. To achieve this last
condition, it suffices to choose a $\mathcal{
C}^{2,\alpha}$-smooth hypersurface $M_1$ passing through $p_1$ with
$T_{p_1}M_0=T_{p_1}M_1$ such that $M_1\backslash \{p_1\}$ is contained
in $M_0^+$. 
\hfill $\square$

\section{Analytic discs}

Let $p_1$ be as in Lemma~2.1. First, we can choose coordinates vanishing at 
$p_1$ and represent $M$ near $p_1$ by the vectorial equation
\def\theequation{3.1}
\begin{equation}
x=h(w,y), \ \ \ \ \ w\in\C^m, \ \ z=x+iy\in\C^n,
\end{equation}
where $h=(h_1,\dots,h_n)$ is of class $\mathcal{ C}^{2,\alpha}$ and
satisfies $h_j(0)=0$ and $dh_j(0)=0$.

Let us recall some generalities ({\it see}~[Bo] for
background).  Denote by $\Delta$ the open unit disc in $\C$.  An
{\it analytic disc} attached to $M$ is a holomorphic mapping
$A:\Delta\rightarrow\C^N$ which extends continuously (or
$\mathcal{ C}^{k,\alpha}$-smoothly) up to the
boundary $\partial\Delta$ and fulfills $A(\partial\Delta)\subset
M$.

Discs of small size (for example with respect to the $\mathcal{
C}^{2,\alpha}$-norm, $0<\alpha<1$) which are attached to $M$
are then obtained as the solutions of the (modified) Bishop equation
\def\theequation{3.2}
\begin{equation}
Y=T_1[h(W(\cdot),Y(\cdot))]+y_0,
\end{equation}
where $T_1$ denotes the harmonic conjugate operator (Hilbert transform
on $\partial\Delta$) normalized at $\zeta=1$, namely satisfying $T_1u(1)=0$
for any $u\in\mathcal{C}^{2,\alpha}(b\Delta,\R^n)$. One verifies that
every small $\mathcal{ C}^{2,\alpha}$-smooth disc
$A(\zeta)=(W(\zeta),Z(\zeta))=(W(\zeta),X(\zeta)+iY(\zeta))$ attached
to $M$ satisfies~(3.2). Conversely, for $W(\zeta)$ of small $\mathcal{
C}^{2,\alpha}$-norm, equation~(3.2) possesses a unique solution
$Y(\zeta)$, and one easily checks that
$A(\zeta):=(W(\zeta),h(W(\zeta),Y(\zeta))+iY(\zeta))$ is then the unique
disc attached to $M$ with $Y(1)=y_0$ and $w$-component equal to
$W(\zeta)$. According to an optimal analysis of the regularity of
Bishop's equation due to Tumanov [Tu2]  (and valid more
generally in the classes $\mathcal{C}^{k,\alpha}$ for $k\geq 1$ and
$0<\alpha<1$), $Y(\zeta)$ and then $A(\zeta)$ are of class
$\mathcal{C}^{2,\alpha}$ over $\overline{\Delta}$.

After a linear transformation we can assume that the tangent space to
$M_1$ is given by $\{x=0, \, u_1=0\}$ and that $T_0M^+$ is given by
$\{u_1>0\}$ near the origin. Let $\rho_0>0$ be small
and let $A$ be the analytic disc we obtain by
solving 
\def\theequation{3.3}
\begin{equation}
Y=T_1 h[(W(\cdot),Y(\cdot))],\ \ \ 
\text{\rm with} \ \ \ W(\zeta):=(\rho_0-\rho_0\zeta,0,
\dots,0).
\end{equation}
Notice that the disc $W_1(\zeta):=(\rho_0-\rho_0\zeta)$ satisfies
$W_1(1)=0$ and $W_1(\overline{\Delta}\backslash \{1\})\subset
\{u_1+iv_1 \in\C :\, u_1>0\}$.  Elementary properties of Bishop's
equation yield $A(\partial\Delta)\backslash \{1\} \subset M^+$ 
if $\rho_0>0$ is
sufficiently small ({\it cf.}  [MP1, Lemma~2.4]). 
{\sc Figure~2} below is devoted to provide a
geometric intuition of the relative situation of the boundary of 
the disc $A$ with respect to $M_1$.

\bigskip
\begin{center}
\begin{picture}(0,0)%
\epsfig{file=figure2.pstex}%
\end{picture}%
\setlength{\unitlength}{3947sp}%
\begingroup\makeatletter\ifx\SetFigFont\undefined%
\gdef\SetFigFont#1#2#3#4#5{%
  \reset@font\fontsize{#1}{#2pt}%
  \fontfamily{#3}\fontseries{#4}\fontshape{#5}%
  \selectfont}%
\fi\endgroup%
\begin{picture}(5874,3759)(334,-3133)
\put(2799,-1498){\makebox(0,0)[lb]{\smash{\SetFigFont{8}{9.6}{\familydefault}{\mddefault}{\updefault}$0=A(1)$}}}
\put(2836,355){\makebox(0,0)[lb]{\smash{\SetFigFont{8}{9.6}{\familydefault}{\mddefault}{\updefault}$w_2,\ldots,w_m,y$}}}
\put(3849,-1130){\makebox(0,0)[lb]{\smash{\SetFigFont{8}{9.6}{\familydefault}{\mddefault}{\updefault}$A(\partial\Delta)$}}}
\put(5738,-1295){\makebox(0,0)[lb]{\smash{\SetFigFont{8}{9.6}{\familydefault}{\mddefault}{\updefault}$u_1$}}}
\put(2987, 32){\makebox(0,0)[lb]{\smash{\SetFigFont{8}{9.6}{\familydefault}{\mddefault}{\updefault}$M_1$}}}
\put(3369,-260){\makebox(0,0)[lb]{\smash{\SetFigFont{8}{9.6}{\familydefault}{\mddefault}{\updefault}$M^+$}}}
\put(2432,-260){\makebox(0,0)[lb]{\smash{\SetFigFont{8}{9.6}{\familydefault}{\mddefault}{\updefault}$M^-$}}}
\put(4080,-455){\makebox(0,0)[lb]{\smash{\SetFigFont{8}{9.6}{\familydefault}{\mddefault}{\updefault}$v_1$}}}
\put(5303,-290){\makebox(0,0)[lb]{\smash{\SetFigFont{8}{9.6}{\familydefault}{\mddefault}{\updefault}$M$}}}
\put(1592,-1434){\makebox(0,0)[lb]{\smash{\SetFigFont{7}{8.4}{\familydefault}{\mddefault}{\updefault}$E_{\rm nr}$}}}
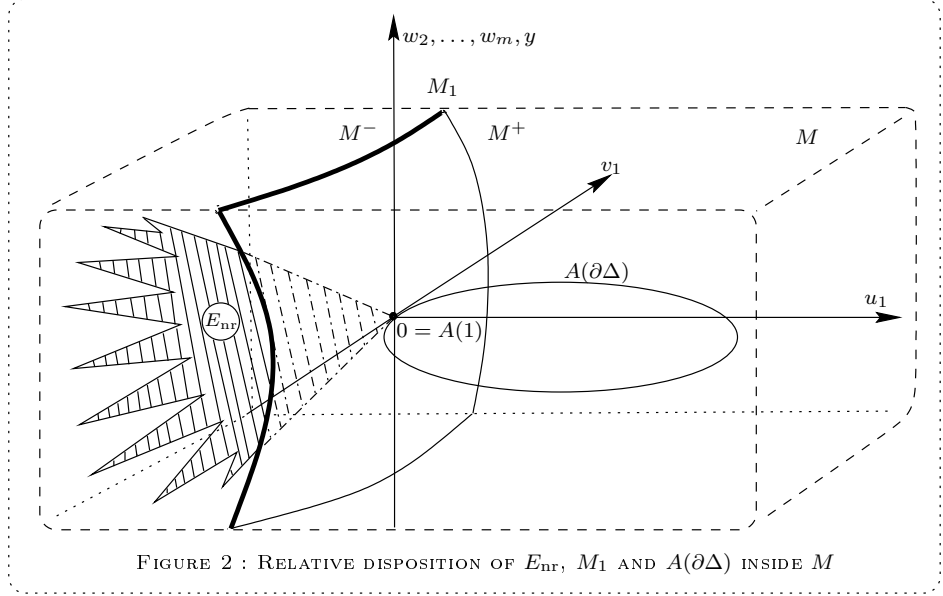
\put(1165,-2956){\makebox(0,0)[lb]{\smash{\SetFigFont{8}{9.6}{\familydefault}{\mddefault}{\updefault}{\sc Figure 2: Relative disposition of $E_{\rm nr}$, $M_1$ and $A(\partial\Delta)$ inside $M$}}}}
\end{picture}

\end{center}
\bigskip

At first, we explain how one usually constructs small wedges attached
to $M$ at $p_1$ by means of {\it deformations of analytic discs} and
then in Sections~4, 5 and 6 below,
we shall explain some of the modifications which are needed in the
presence of a singularity $E_{\rm nr}$ in order to produce wedge extension at
$p_1$.  Following [MP3, pp.~863--864], we shall include (or say, 
``deform'') $A$ in a para\-metrized family $A_{\rho,s,v}$ with
varying radius $\rho$ plus supplementary parameters $s,v$ and with
$A_{\rho_0,0,0}=A$.  During the construction, we shall sometimes
permit ourselves to decrease parameters, related constants,
neighborhoods and domains of existence without explicit mentioning.
At present, our goal is to explain how we can add some conveninent extra 
simplifying assumptions to the hypotheses of Theorem~1.1, 
{\it see} especially conditions {\bf 1)}, {\bf 2)} and
{\bf 3)} before Theorem~3.1 below.

Let $\mathcal{ W}_1$ be the wedgelike domain attached to $M\backslash
E_{\rm nr}$ constructed in Section~2 and let $f\in\mathcal{
O}(\mathcal{ W}_1)$. We want to extend $f$ holomorphically to a wedge
of edge a small neighborhood of the special point $p_1\in E_{\rm nr}$
picked thanks to Lemma~2.1.  Let $\mathcal{ W}_2\subset\mathcal{ W}_1$
be a small wedge attached to a neighborhood of $A(-1)$ in $M^+$. As in
[Tu1,2], [MP1,3], we can construct analytic discs
$A_{\rho,s,v}=(W_{\rho,s,v},Z_{\rho,s,v})$ attached to $M\cup
\mathcal{ W}_2$ with the following properties:

\smallskip
\begin{itemize}
\item[{\bf (1)}]
The parameters $s,v$ belong to neighborhoods $U_s,U_v$ of $0$
in $\R^{2m+n-1}$, $\R^{n-1}$ respectively and $\rho$ belongs to 
the interval $[0,\rho_1)$, for some $\rho_1>\rho_0$.
\item[{\bf (2)}]
The mapping $(\rho,s,v)\mapsto A_{\rho,s,v}$ is of class $\mathcal{
C}^{2,\beta}$ for all $0<\beta<\alpha$. For $\rho\neq 0$, these maps
are embeddings of $\overline{\Delta}$ into $\C^{m+n}$.  Finally, we
have $A_{\rho_0,0,0}=A$ and the discs $A_{0,s,v}$ are constant.
\item[{\bf (3)}]
For every fixed $v_0\in U_v$, the union $\bigcup_{s\in U_s}
A_{\rho_0,s,v_0}(\{e^{i\theta}\, :\, \vert\theta\vert< \pi/4\})$ is
an open subset of $M$ containing the origin which is 
$\mathcal{ C}^{2,\beta}$-smoothly
foliated by the curves $A_{\rho_0,s,v_0}(\{e^{i\theta}\, :\,
\vert\theta\vert< \pi/4\})$.
\item[{\bf (4)}]
The mapping $U_v\ni v\mapsto [{d\over d\theta}
A_{\rho_0,0,v}(e^{i\theta})]_{\theta=0}\in T_0M/T_0^cM\simeq \R^n$ has
rank $n-1$ and its image is transverse to the vector $[{d\over
d\theta} A(e^{i\theta})]_{\theta=0}\in T_0M/T_0^cM\simeq \R^n$.
In geometric terms, this property means that the union of 
tangent real lines
$$
\R\,\left[{d\over d r}
A_{\rho_0,0,v}(r e^{i\theta})\right]_{\zeta=1}=
-i\,\R\,\left[{d\over d\theta}
A_{\rho_0,0,v}(e^{i\theta})\right]_{\theta=0}
$$
spans an open cone in the normal bundle to $M$, namely 
$T_0\C^{m+n}/T_0M\cong i(T_0M/T_0^cM)$.
\item[{\bf (5)}]
Let $\omega=\{\zeta\in\Delta :\,|\zeta-1|<\delta\}$ be a neighborhood
of $1$ in $\Delta$, with some small $\delta>0$. It follows from
properties {\bf (3)} and {\bf (4)} that the union $\mathcal{
W}=\bigcup_{s\in U_s,v\in U_v} A_{\rho_0,s,v}(\omega)$ is an open
wedge of edge a neighborhood of the origin in $M$ which is foliated by
the discs $A_{\rho_0,s,v}(\omega)$.
\item[{\bf (6)}]
The sets $D_{s,v}=\bigcup_{0\leq\rho<\rho_1,|\zeta|=1}
A_{\rho,s,v}(\zeta)$ are {\it real} two-dimensional discs of class
$\mathcal{ C}^{2,\beta}$ embedded in $M$ which are foliated (with a
circle degenerating to a point for $\rho=0$) by the circles
$A_{\rho,s,v}(\partial\Delta)$.
\item[{\bf (7)}]
There exists a $(2m+n-2)$-dimensional submanifold $H$ of $\R^{2m+n-1}$
passing through the origin such that for every fixed $v_0\in U_v$, the union
$\bigcup_{s\in H} \, D_{s,v_0}$ is a $(\dim\, M)$-dimensional open box
foliated by real $2$-discs which is contained in $M$ and which contains
the origin. Intruitively, it is a stack of plates.
\end{itemize}

\smallskip
Let us make some commentaries.  We stress that the family
$A_{\rho,s,v}$ is obtained by solving the Bishop equation for
explicitly prescribed data ({\it see} [MP3, p.~837] or [MP1, p.~863];
the important Lemma~2.7 in [MP1] which produces the parameter $v$
satisfying {\bf (4)} above is due to Tumanov [Tu1]). Since
Bishop's equation is very flexible, this entails that
every geometrical property of the family is stable under slight
perturbation of the data.  Notice for instance that as $A$ is an
embedding of $\overline{\Delta}$ into $\C^{m+n}$, all its small
deformations will stay embeddings.  In particular we get a likewise
family $A^d_{\rho,s,v}$ if we replace $M$ by a slightly deformed
$\mathcal{ C}^{2,\alpha}$-smooth manifold $M^d$ (this corresponds to
replacing $h$ by a function $h^d$ close to $h$ in $\mathcal{
C}^{2,\alpha}$-norm in~(3.1), (3.2) and~(3.3)).

\smallskip
\noindent
{\it Further remark.}
If $A'$ is an arbitrary disc which is sufficiently close to $A$
in $\mathcal{C}^{1,\beta}$-norm, for some $0<\beta<\alpha$, we can
also include $A'$ in a similar $\mathcal{C}^{1,\gamma}$-smooth
($0<\gamma<\beta$) family $A_{\rho,s}'$, {\it without the parameter $v$},
which satisfies the geometric properties {\bf (3)}, {\bf (6)}
and {\bf (7)}
above. This remark will be useful in the end of Section~4 below.

\smallskip
Using such a nice family $A_{\rho,s,v}$ which gently deforms as a
family $A_{\rho,s,v}^d$ under perturbations, let us begin to remind
from [MP1] how we can add three simplifying geometric assumptions to
Theorem~1.1, without loss of generality.

First of all, using apartition of unity, we can perform arbitrarily
small $\mathcal{ C}^{2,\alpha}$-smooth deformations $M^d$ of $M$
leaving $E_{\text{\rm nr}}$ fixed and moving $M\backslash E_{\text{\rm
nr}}$ inside the wedgelike domain $\mathcal{ W}_1$.  Further, we can
make $M^d$ to depend on a single small real parameter $d\geq 0$ with
$M^0=M$ and $M^d\backslash E_{\rm nr}\subset \mathcal{ W}_1$ for all
$d>0$. Now, {\it the wedgelike domain $\mathcal{ W}_1$ becomes a
neighborhood of $M^d$ in $\C^{m+n}$}.  In the sequel, we shall denote
this neighborhood by $\Omega$. By stability of Bishop's equation, we
obtain a deformed disc $A^d$ attached to $M^d$ by solving~(3.3) with
$h^d$ in place of $h$.  In the sequel, we will also consider a small
neighborhood $\Omega_1$ of $A^d(-1)$ in $\C^{m+n}$ which contains the
intersection of the above wedge $\mathcal{ W}_2$ with a neighborhood
of $A(-1)$ in $\C^{m+n}$.

Again by stability of Bishop's equation, we also obtain deformed
families $A^d_{\rho,s,v}$ attached to $M^d\cup \Omega_1$, satisfying
properties {\bf (1)-(7)} above. Recall that according to [Tu2], the
mapping $(\rho,s,v,d)\mapsto A^d_{\rho,s,v}$ is $\mathcal{
C}^{2,\beta}$-smooth for all $0<\beta<\alpha$. In the core of the proof of
our main Theorem~1.1 (Sections 4, 5 and 6 below), we will show that,
for each sufficiently small fixed $d>0$, we get holomorphic extension to
the wedgelike set $\mathcal{ W}^d=\bigcup_{s\in U_s,v\in U_v}
A^d_{\rho_0,s,v}(\omega)$ attached to a neighborhood of $0$ in
$M^d$. But this implies Theorem~\ref{main}: In the limit $d\rightarrow
0$, the wedges $\mathcal{ W}^d$ tend smoothly to the wedge $\mathcal{
W}:=\mathcal{ W}^0$ attached to a neighborhood of $0$ in $M^0=M$.  As
the construction depends smoothly on the deformations $d$, we derive
univalent holomorphic extension to $\mathcal{ W}$ thereby arriving at
a contradiction to the definition of $E_{\text{\rm nr}}$.

As a summary of the above discussion, we formulate below the local
statement that remains to prove. Essentially, we have shown that it
suffices to prove Theorem~1.1 with the following three extra simplifying
assumptions: 
\begin{itemize}
\item[{\bf 1)}] 
Instead of functions which are holomorphic in a
wedgelike open set attached to $M\backslash E_{\rm nr}$, we consider
functions which are holomorphic in a neighborhood of $M\backslash
E_{\rm nr}$ in $\C^{m+n}$. 
\item[{\bf 2)}]
Proceeding by contradiction, we
have argued that it suffices to remove at least one point of $E_{\rm
nr}$. 
\item[{\bf 3)}]
Moreover, we can assume that the point $p_1\in E_{\rm
nr}$ we want to remove is behind a generic ``wall'' $M_1$ as depicted
in {\sc Figure~2}.
\end{itemize}

Consequently, from now on, we shall denote the set $E_{\rm nr}$ simply
by $E$. We also denote $M^d$ simply by $M$. We take again the disc $A$
defined by~(3.3) and its deformation $A_{\rho,s,v}$. The goal is now
to show that holomorphic functions in a neighborhood of $M\backslash
E$ in $\C^{m+n}$ extend holomorphically to a wedge at $p_1$, assuming
the ``nice'' geometric situation of {\sc Figure~2}. To be precise, we have
argued that Theorem~1.1 is reduced to the following precise
and geometrically more concrete statement.
 
\begin{thm}\label{red}
Let $M$ be a $\mathcal{C}^{2,\alpha}$-smooth generic CR manifold in
$\C^{m+n}$ of codimension $n$. Let $M_1\subset M$ be a
$\mathcal{C}^{2,\alpha}$-smooth generic CR manifold of dimension
$2m+n-1$ and let $p_1\in M_1$. Let $M^+$ and $M^-$ denote the two
local open sets in which $M$ is divided by $M_1$, in a neighborhood of
$p_1$. Suppose that $E\subset M$ is a \text{\rm nonempty} closed
subset with $p_1\in E$ satisfying the Hausdorff condition
$H^{2m+n-2}(E)=0$ and suppose that $E\subset M^-\cup\{p_1\}$ {\sc
({\sc Figure~2})}. Let $\Omega$ be a neighborhood of $M\backslash E$ in
$\C^{m+n}$, let $A$ be the disc defined by~(3.3), let $\Omega_1$ be a
neighborhood of $A(-1)$ in $\C^{m+n}$ which is contained in $\Omega$
and let $A_{\rho,s,v}$ be a family of discs attached to $M\cup
\Omega_1$ with the properties {\bf (1)-(7)} explained above. Then
every function $f$ which is holomorphic in $\Omega$ extends
holomorphically to the wedge $\mathcal{ W}=\bigcup_{s\in U_s,v\in U_v}
A_{\rho_0,s,v}(\omega)$.
\end{thm}

Of course, Theorem~3.1 would be obvious if $E$ would be
empty, but we have to take account of $E$.

\section{Proof of Theorem~\ref{red}, part I}

This section contains the part of the proof of Theorem~\ref{red} above
which relies on constructions with the small discs $A_{\rho,s,v}$
attached to $M\cup \Omega$. Since we want the boundaries of our discs
to avoid $E$, we shall employ the following elementary lemma several
times, which is simply a convenient particularization of
a general property of Hausdorff measures [C, Appendix A6].

\begin{lem}
Let $N$ be a real $d$-dimensional manifold and let $E\subset N$ be a
closed subset. Let $U$ be a small neighborhood of
the origin in $\R^{d-1}$ and let $\Phi:\partial\Delta\times U\to N$
$($resp. $\Psi:(0,1)\times U\to M)$ be an embedding. 
\begin{itemize}
\item[{\bf (i)}]
If $H^{d-2}(E)=0$, then the set of $x\in U$ for which
$\Phi(\partial\Delta\times\{x\})\cap E$ is nonempty
$($resp. $\Psi((0,1)\times\{x\})\cap E\neq\emptyset)$ is of
zero $(d-2)$-dimensional Hausdorff measure.
\item[{\bf (ii)}] 
If $H^{d-1}(E)=0$, then for almost every $x\in U$ in the sense of
Lebesgue measure, we have $\Phi(\partial\Delta\times\{x\})\cap
E=\emptyset$ $($resp. $\Psi((0,1)\times\{x\})\cap
E=\emptyset)$.
\end{itemize}
\end{lem}

\smallskip
\noindent    
{\it Proof of Theorem~\ref{red}:}  
We divide the proof in five steps.                                      
     
\smallskip
\noindent    
{\bf Step 1: Holomorphic extension to a dense subset of $\mathcal{ W}$.}   
We shall start by constructing a holomorphic extension to an
everywhere dense open subdomain
of the wedge $\mathcal{ W}=\bigcup_{s\in U_s,v\in U_v}
A_{\rho_0,s,v}(\omega)$ by means of the disc technique
(continuity principle).

For each fixed $v_0\in U_v$, the first
dimensional count of Lemma~4.1 (which applies by the foliation property
{\bf (3)} of the discs) yields a closed subset $\mathcal{
S}_{v_0}\subset U_s$ depending on $v_0$ and satisfying
$H^{2m+n-2}(\mathcal{ S}_{v_0})=0$ such that for every
$s\notin\mathcal{ S}_{v_0}$ we have $A_{\rho_0,s,v_0}(\partial\Delta)\cap
E=\emptyset$. Notice also that $\mathcal{ S}_{v_0}$ does not locally
disconnect $U_s$, for dimensional reasons ([C, Appendix A6]).

By property~{\bf (7)} of Section~3, the real two-dimensional 
discs $D_{s,v_0}$ foliate an open subset of $M$, for $s$ 
running in a manifold $H$ of dimension $2m+n-2$. Consequently, for almost
every $s\in H$, (in the sense of Lebesgue measure), we have
$D_{s,v_0}\cap E=\emptyset$.

\bigskip
\hspace{-1.4cm}
\begin{picture}(0,0)%
\epsfig{file=figure3.pstex}%
\end{picture}%
\setlength{\unitlength}{3947sp}%
\begingroup\makeatletter\ifx\SetFigFont\undefined%
\gdef\SetFigFont#1#2#3#4#5{%
  \reset@font\fontsize{#1}{#2pt}%
  \fontfamily{#3}\fontseries{#4}\fontshape{#5}%
  \selectfont}%
\fi\endgroup%
\begin{picture}(6579,3069)(184,-2338)
\put(5857,-1057){\makebox(0,0)[lb]{\smash{\SetFigFont{6}{7.2}{\familydefault}{\mddefault}{\updefault}$\simeq$}}}
\put(5310,450){\makebox(0,0)[lb]{\smash{\SetFigFont{6}{7.2}{\familydefault}{\mddefault}{\updefault}$\partial\Delta$}}}
\put(6175, 14){\makebox(0,0)[lb]{\smash{\SetFigFont{6}{7.2}{\familydefault}{\mddefault}{\updefault}$\omega$}}}
\put(5282,  1){\makebox(0,0)[lb]{\smash{\SetFigFont{6}{7.2}{\familydefault}{\mddefault}{\updefault}$-1$}}}
\put(5886,-510){\makebox(0,0)[lb]{\smash{\SetFigFont{6}{7.2}{\familydefault}{\mddefault}{\updefault}$-i$}}}
\put(5906,540){\makebox(0,0)[lb]{\smash{\SetFigFont{6}{7.2}{\familydefault}{\mddefault}{\updefault}$i$}}}
\put(6469,  5){\makebox(0,0)[lb]{\smash{\SetFigFont{6}{7.2}{\familydefault}{\mddefault}{\updefault}$1$}}}
\put(5682,238){\makebox(0,0)[lb]{\smash{\SetFigFont{6}{7.2}{\familydefault}{\mddefault}{\updefault}$\Delta$}}}
\put(5490,-1799){\makebox(0,0)[lb]{\smash{\SetFigFont{6}{7.2}{\familydefault}{\mddefault}{\updefault}$\partial\Delta$}}}
\put(5816,-1972){\makebox(0,0)[lb]{\smash{\SetFigFont{6}{7.2}{\familydefault}{\mddefault}{\updefault}$\Delta$}}}
\put(6260,-1951){\makebox(0,0)[lb]{\smash{\SetFigFont{6}{7.2}{\familydefault}{\mddefault}{\updefault}$\omega$}}}
\put(6466,-1796){\makebox(0,0)[lb]{\smash{\SetFigFont{6}{7.2}{\familydefault}{\mddefault}{\updefault}$\partial\Delta$}}}
\put(4955,-763){\makebox(0,0)[lb]{\smash{\SetFigFont{6}{7.2}{\familydefault}{\mddefault}{\updefault}$A_{\rho_0,s,v}$}}}
\put(3408,410){\makebox(0,0)[lb]{\smash{\SetFigFont{6}{7.2}{\familydefault}{\mddefault}{\updefault}$A_{\rho_0,s,v}(\Delta)$}}}
\put(5130,-1604){\makebox(0,0)[lb]{\smash{\SetFigFont{6}{7.2}{\familydefault}{\mddefault}{\updefault}$M$}}}
\put(639,-1788){\makebox(0,0)[lb]{\smash{\SetFigFont{6}{7.2}{\familydefault}{\mddefault}{\updefault}$\Omega$}}}
\put(1624,-457){\makebox(0,0)[lb]{\smash{\SetFigFont{6}{7.2}{\familydefault}{\mddefault}{\updefault}$E_{\mathcal{ W}}$}}}
\put(973,222){\makebox(0,0)[lb]{\smash{\SetFigFont{6}{7.2}{\familydefault}{\mddefault}{\updefault}$\C^{m+n}$}}}
\put(1999,-1302){\makebox(0,0)[lb]{\smash{\SetFigFont{6}{7.2}{\familydefault}{\mddefault}{\updefault}$M_1$}}}
\put(1794,-1424){\makebox(0,0)[lb]{\smash{\SetFigFont{6}{7.2}{\familydefault}{\mddefault}{\updefault}$A_{\rho_0,s,v}(\partial\Delta)$}}}
\put(1531,-1261){\makebox(0,0)[lb]{\smash{\SetFigFont{6}{7.2}{\familydefault}{\mddefault}{\updefault}$E$}}}
\put(981,-1467){\makebox(0,0)[lb]{\smash{\SetFigFont{6}{7.2}{\familydefault}{\mddefault}{\updefault}$E$}}}
\put(1198,-2041){\makebox(0,0)[lb]{\smash{\SetFigFont{6}{7.2}{\familydefault}{\mddefault}{\updefault}{\sc Figure 3: Description of $E_{\mathcal{ W}}:=$ union of analytic discs }}}}
\put(1194,-2152){\makebox(0,0)[lb]{\smash{\SetFigFont{6}{7.2}{\familydefault}{\mddefault}{\updefault}{\sc $A_{\rho_0,s,v}(\Delta)$ whose boundary $A_{\rho_0,s,v}(\partial\Delta)$ meets the singularity $E$}}}}
\put(2582,-136){\makebox(0,0)[lb]{\smash{\SetFigFont{6}{7.2}{\familydefault}{\mddefault}{\updefault}$\mathcal{ W}$}}}
\put(4565,-1573){\makebox(0,0)[lb]{\smash{\SetFigFont{6}{7.2}{\familydefault}{\mddefault}{\updefault}$\Omega$}}}
\put(308,-1789){\makebox(0,0)[lb]{\smash{\SetFigFont{6}{7.2}{\familydefault}{\mddefault}{\updefault}$M$}}}
\put(4405,-1216){\makebox(0,0)[lb]{\smash{\SetFigFont{6}{7.2}{\familydefault}{\mddefault}{\updefault}$A_{\rho_0,s,v}(\partial\Delta)$}}}
\put(2337,-1212){\makebox(0,0)[lb]{\smash{\SetFigFont{6}{7.2}{\familydefault}{\mddefault}{\updefault}$M^+$}}}
\put(1539,-1428){\makebox(0,0)[lb]{\smash{\SetFigFont{6}{7.2}{\familydefault}{\mddefault}{\updefault}$M^-$}}}
\put(1871,-927){\makebox(0,0)[lb]{\smash{\SetFigFont{6}{7.2}{\familydefault}{\mddefault}{\updefault}$A_{\rho_0,s,v}(\omega)$}}}
\end{picture}

\bigskip

Since $E$ is closed, we claim that for every $s\notin\mathcal{
S}_{v_0}$, it follows that we can contract every boundary
$A_{\rho_0,s,v_0}(\partial \Delta)$ which does not
meet $E$, to a point in $M$ without meeting
$E$ by an analytic isotopy ({\it cf.}~[MP3, p.~864]).  Indeed, by
shifting $s$ to some nearby $s'$, we
first move $A_{\rho_0,s,v_0}$ into a disc
$A_{\rho_0,s',v_0}$ which also satisfies $A_{\rho_0,s',v_0}(\partial
\Delta)\cap E=\emptyset$. Choosing well $s'$, this boundary belongs to 
a real disc $D_{s',v_0}$ satisfying
$D_{s',v_0}\cap E=\emptyset$. This can be achieved with $s'$
arbitrarily close to $s$, since $\mathcal{ S}_{v_0}$ does not
disconnect $U_s$. Then we contract in the obvious manner the disc
$A_{\rho_0,s',v_0}$ to the point $A_{0,s',v_0}(\overline{\Delta})$ by
isotoping its boundary inside $D_{s',v_0}$ (recall that 
$D_{s',v_0}$ is a union of boundary of discs).  Applying the continuity
principle to this analytic isotopy of discs, we see that we can extend
every function $f\in\mathcal{O}(\Omega)$ holomorphically
to a neighborhood of $A_{\rho_0,s,v_0}(\overline{\Delta})$ 
in $\C^{m+n}$,
for every $s\not\in \mathcal{ S}_{v_0}$ and for every $v_0\in U_v$.

From the nice geometry {\bf (5)} of the family $A_{\rho,s,v}$ one
easily derives that the various local extensions near
$A_{\rho_0,s,v_0}(\omega)$ for $s\notin\mathcal{ S}_{v_0}$ fit in a
univalent function $F\in\mathcal{O}(\mathcal{ W}\backslash E_\mathcal{
W})$, where $E_\mathcal{ W}:=\bigcup_{s\in\mathcal{ S}_{v_0}, v_0\in
U_v} A_{\rho_0,s,v_0}(\omega)$.  Furthermore we observe that $E_\mathcal{
W}$ is laminated by holomorphic discs and satisfies
$H^{2m+2n-1}(E_\mathcal{ W})=0$. This metrical property implies that
$\mathcal{ W}\backslash E_\mathcal{ W}$ is locally connected.  The remainder
of the proof is devoted to show how to extend $F$ through
$E_{\mathcal{W}}$. This occupies the paper up to its end.  The
difficulty and the length of the proof comes from the fact that the
disc method necessarily increases by a factor $1$ the dimension of the
singularity: it transforms a singularity set $E\subset M$ of
codimension $2^{+0}$ into a bigger singularity set $E_{\mathcal{
W}}\subset \mathcal{ W}$ which is of codimension $1^{+0}$.

\smallskip
\noindent
{\bf Step 2: Plan for the removal of $E_{\mathcal{W}}$.} Let us
remember that our goal is to show that $p_1$ is ${\mathcal
W}$-removable in order to achieve the final step in our reasoning by
contradiction which begins in Section~2.  To show that $p_1$ is
removable, it suffices to extend $F$ through $E_{\mathcal{ W}}$. At
first, we notice that because $H^{2m+2n-1}(E_{\mathcal{ W}})=0$, it
follows that $\mathcal{ W}\backslash E_{\mathcal{ W}}$ is locally
connected, so the part of the envelope of holomorphy of $\mathcal{
W}\backslash E_{\mathcal{ W}}$ which is contained in $\mathcal{ W}$ is
not multisheeted: it is necessarily a subdomain of $\mathcal{ W}$. In
analogy with the beginning of Section~2, let us therefore denote by
$E_\mathcal{ W}^{\rm nr}$ the set of points of $E_{\mathcal{ W}}$
through which our holomorphic function $F\in\mathcal{ O}(\mathcal{
W}\backslash E_{\mathcal{ W}})$ does not extend holomorphically.  If
$E_\mathcal{ W}^{\rm nr}$ is empty, we are done, gratuitously.  As it
might certainly be nonempty, we shall suppose therefore that
$E_\mathcal{ W}^{\rm nr}\neq \emptyset$ and we shall construct a
contradiction in the remainder of the paper. Let $q\in E_\mathcal{
W}^{\rm nr}\neq \emptyset$. To derive a contradiction, it suffices to
show that $F$ extends holomorphically through $q$.  Philosophically
again, it will suffice to remove one single point, which will simplify
the presentation and the geometric reasonings. Finally, as
$E_\mathcal{ W}^{\rm nr}\neq \emptyset$ is contained in $E_{\mathcal{
W}}$, there exist a point $\zeta_0\in\partial\Delta$ and parameters
$(\rho_0,s_0,v_0)$ such that $q=A_{\rho_0,s_0,v_0}(\zeta_0)$.  In the
sequel, we shall simply denote the disc $A_{\rho_0,s_0,v_0}$ by
$A_{\rm nr}$. Obviously also, $H^{2m+2n-1}(E_{\mathcal{ W}}^{\rm nr})=0$.

\smallskip
\noindent
{\bf Step 3: Smoothing the boundary of the singular disc $A_{\rm nr}$ 
near $\zeta=-1$.}  In step~4 below, our goal will be to deform $A_{\rm
nr}$ to extend $F$ through $q$. As we shall need to glue a maximally
real submanifold $R_1$ of $M$ along $A_{\rm
nr}(\partial\Delta\backslash \{\v\zeta+1\v<\eps\})$ to some collection
of maximally real {\it planes} along $A_{\rm nr}(\zeta)$ for
$\zeta\in\partial\Delta$ near $-1$, and because
$\mathcal{C}^{2,\beta}$-smoothness of $A_{\rm nr}$ will not be
sufficient to keep the $\mathcal{C}^{2,\beta}$-smoothness of the glued
object, it is convenient to smooth out first $A_{\rm nr}$ near
$\zeta=-1$ ({\it see}~especially Step~2 of Section~6 below).
Fortunately, we can use the freedom $\Omega_1$ (the small
neighborhood of $A(-1)$ in Theorem~3.1) to modify the boundary
of $A_{\rm nr}$. Thus, {\it for technical reasons only}, we need the
following preliminary lemma, which is simply obtained by
reparametrizing an almost full subdisc of $A_{\rm nr}$.  This
preparatory reparametrization is indispensible to state our Main
Lemma~4.3 below correctly.

\begin{lem}
For every $\eps>0$, there exists an analytic disc $A'$ satisfying
\begin{itemize}
\item[{\bf (a)}]
$A'$ is a $\mathcal{ C}^{2,\beta}$-smooth subdisc of $A_{\rm nr}$, namely 
$A'(\overline{\Delta})\subset A_{\rm nr}(\overline{\Delta})$,
such that 
moreover $A'(\overline{\Delta})\supset
A_{\rm nr}(\overline{\Delta}
\backslash \{\v\zeta+1\v<2\eps\})$.
\item[{\bf (b)}]
$A'$ is real analytic over $\{\zeta\in\partial\Delta: 
\v \zeta+1\v<\eps\}$.
\item[{\bf (c)}]
$\v\v A'-A_{\rm nr}\v\v_{\mathcal{ C}^{2,\beta}}\leq \eps$.
\item[{\bf (d)}]
$A'(\partial\Delta)\subset M\cup \Omega_1$.
\end{itemize}
\end{lem}

\proof
Of course, {\bf (d)} follows immediately from {\bf (a)} and {\bf (c)}
if $\eps$ is sufficiently small. To construct $A'$, we consider a
$\mathcal{C^\infty}$-smooth cut-off function
$\mu_\eps:\partial\Delta\to [0,1]$ with $\mu_\eps(\zeta)=1$ for
$\v\zeta+1\v>2\eps$ and $\mu_\eps(\zeta)$ equal to a constant
$c_\eps<1$ with $c_\eps >1-\eps$ for $\v\zeta+1\v<\eps$. Let
$\Delta_{\mu_\eps}$ be the (almost full) subdisc of $\Delta$ defined
by $\{\zeta\in\Delta: \v\zeta\v<\mu_\eps(\zeta/\v\zeta\v)\}$. Let
$\psi_\eps$ be the Riemann conformal map $\Delta\to
\Delta_{\mu_\eps}$. We can assume that $\psi_\eps(-1)=-c_\eps\in
\partial\Delta_{\mu_\eps}\cap \R$.
By Caratheodory's theorem and by the Schwarz
symmetry principle, $\psi_\eps$ is $\mathcal{ C}^\infty$-smooth up to the
boundary and real analytic near $\zeta=-1$.  If $\eps$ is sufficiently
small and $c_\eps$ sufficiently close to $1$, the stability of
Riemann's uniformization theorem under small $\mathcal{C}^\infty$-smooth
perturbations shows that the disc
$$
A'(\zeta):=A_{\rm nr}(\psi_\eps(\zeta)).
$$
satisfies the desired properties, possibly with a slightly different
small $\eps$.
\endproof

\smallskip
\noindent
{\bf Step 4: Variation of the singular disc.} 
In the sequel, we shall constantly denote the disc of Lemma~4.2 by
$A'$.  We set $\zeta_q:=\psi_\eps^{-1}(\zeta_0)$, so that
$A'(\zeta_q)=q$. Of course, after a reparametrization by a Blaschke
transformation, we can (and we will) assume that $\zeta_q=0$.  By construction,
$A'\v_{\partial\Delta}$ is real analytic near $-1$ and the point
$q=A'(0)$ is contained $E_{\mathcal W}^{\rm nr}$, the set through
which our partial extension $F$ does not extend {\it a priori}. To
derive a contradiction, our next purpose is to produce a disc $A''$
close to $A'$ and {\it passing through the fixed point $q$} such that
$q$ can be encircled by a small closed curve in $A''(\Delta)\backslash
E_{\mathcal{ W}}^{\rm nr}$, because in such a situation, we will be able to
apply the continuity principle as in the typical local situation of
Hartog's theorem ({\it see} {\bf (4)} of Lemma~4.3 and Step~5 below).

At first glance it seems that we can produce $A''$ simply by turning
$A'$ a little around $q$: indeed, Lemma~4.1 applies, since
$H^{2m+2n-1}(E_{\mathcal{ W}}^{\rm nr})=0$. However, the difficult
point is to guarantee that $A''$ is still attached to the union of $M$
with the small neighborhood $\Omega_1$ of $A(-1)$ in $\C^{m+n}$. The
following key lemma asserts that these additional requirements can be
fulfilled.

\begin{mainlem}\label{tr}
Let $A'$ be the disc of Lemma~4.1, let $q=A'(0)\in
E_{\mathcal{ W}}^{\rm nr}$ and let $0<\beta<\alpha$ be 
arbitrarily close to $\alpha$. 
Then there exists a parameterized family $A_{t'}'$ of
analytic discs with the following properties:
\begin{itemize}
\item[{\bf (1)}]
The parameter $t'$ ranges in a neighborhood $U_{t'}$ of $0$ in
$\R^{2m+2n-1}$ and $A_0'=A'$.
\item[{\bf (2)}]
The mapping $U_{t'}\times\overline{\Delta}\ni(t',\zeta) \mapsto
A_{t'}'(\zeta)\in\C^{m+n}$ is of class $\mathcal{ C}^{1,\beta}$
and each $A_{t'}'$ is an embedding of $\overline{\Delta}$ into 
$\C^{m+n}$.
\item[{\bf (3)}] 
For all $t'\in U_{t'}$, the point $q=A_{t'}'(0)$ is fixed and
$A_{t'}'(\partial\Delta)\subset M\cup\Omega_1$.
Furthermore, there exists a small $\delta>0$ such that the
large boundary part $A_{t'}'(\partial\Delta
\backslash \{\v\zeta+1\v<\delta\})$ is attached to a fixed maximally real 
$(m+n)$-dimensional
$\mathcal{ C}^{2,\alpha}$-smooth submanifold $R_1$ of $M$.
\item[{\bf (4)}] 
For every fixed $\rho_\eps >0$ which is sufficiently small and for $t'$
ranging in a sufficiently small neighborhood of the origin, the union
of circles
$$
\bigcup_{t'} \, \{A_{t'}'(\rho_\eps e^{i\theta}): \, \theta\in\R\}
$$
foliates a neighborhood in $\C^{m+n}$ of the small fixed circle
$\{A'(\rho_\eps e^{i\theta}):\, \theta\in\R\}$ which encircles the
point $q$ inside $A'(\Delta)$. Consequently, by Lemma~4.1, 
for almost all $t'\in U_{t'}$, the circle 
$\{A_{t'}'(\rho_\eps e^{i\theta}): \, \theta\in\R\}$
does not meet $E_{\mathcal{ W}}^{\rm nr}$.
\end{itemize}
\end{mainlem}

Let us make some explanatory commentaries.  Notice that the discs are only
$\mathcal{ C}^{1,\beta}$-smooth, because the underlying method of
Sections~5 and~6 (implicit function theorem in Banach spaces, {\it
cf.}~[G1]) imposes a real loss of smoothness. If we could have produce
a $\mathcal{ C}^{2,\beta}$-smooth family (assuming for instance that
$M$ was $\mathcal{ C}^{3,\alpha}$-smooth from the beginning, or
asking whether the regularity methods of [Tu2] are 
applicable to the {\it global} Bishop equation), we would
have constructed a slightly different family 
and stated instead of {\bf (4)} the following conic-like differential
geometric property:
\begin{itemize}
\item[{\bf (4')}] 
The parameter $t'$ ranges over a neighborhood $U_{t'}$ of the origin in
$\R^{2m+2n-2}$ with $A_{t'}'(0)=q$ for all $t'$ and the mapping 
$$
U_{t'}\ni t' \ \mapsto \ [\partial A_{t'}'/\partial\zeta](0)\in
T_q\C^{m+n}
$$ 
has maximal rank at $t'=0$ with its image being transverse to
the tangent space of $A'(\Delta)$ at $q$.
\end{itemize}

\bigskip
\hspace{-0.7cm}
\begin{picture}(0,0)%
\epsfig{file=figure4.pstex}%
\end{picture}%
\setlength{\unitlength}{3947sp}%
\begingroup\makeatletter\ifx\SetFigFont\undefined%
\gdef\SetFigFont#1#2#3#4#5{%
  \reset@font\fontsize{#1}{#2pt}%
  \fontfamily{#3}\fontseries{#4}\fontshape{#5}%
  \selectfont}%
\fi\endgroup%
\begin{picture}(5979,3834)(102,-3088)
\put(5837,-2262){\makebox(0,0)[lb]{\smash{\SetFigFont{8}{9.6}{\familydefault}{\mddefault}{\updefault}$M$}}}
\put(1032,-153){\makebox(0,0)[lb]{\smash{\SetFigFont{8}{9.6}{\familydefault}{\mddefault}{\updefault}$\C^{m+n}$}}}
\put(2218,-1913){\makebox(0,0)[lb]{\smash{\SetFigFont{8}{9.6}{\familydefault}{\mddefault}{\updefault}$M_1$}}}
\put(1676,-1865){\makebox(0,0)[lb]{\smash{\SetFigFont{8}{9.6}{\familydefault}{\mddefault}{\updefault}$E$}}}
\put(5184,-2226){\makebox(0,0)[lb]{\smash{\SetFigFont{8}{9.6}{\familydefault}{\mddefault}{\updefault}$\Omega$}}}
\put(263,-2475){\makebox(0,0)[lb]{\smash{\SetFigFont{8}{9.6}{\familydefault}{\mddefault}{\updefault}$M$}}}
\put(2608,-1809){\makebox(0,0)[lb]{\smash{\SetFigFont{8}{9.6}{\familydefault}{\mddefault}{\updefault}$M^+$}}}
\put(1686,-2059){\makebox(0,0)[lb]{\smash{\SetFigFont{8}{9.6}{\familydefault}{\mddefault}{\updefault}$M^-$}}}
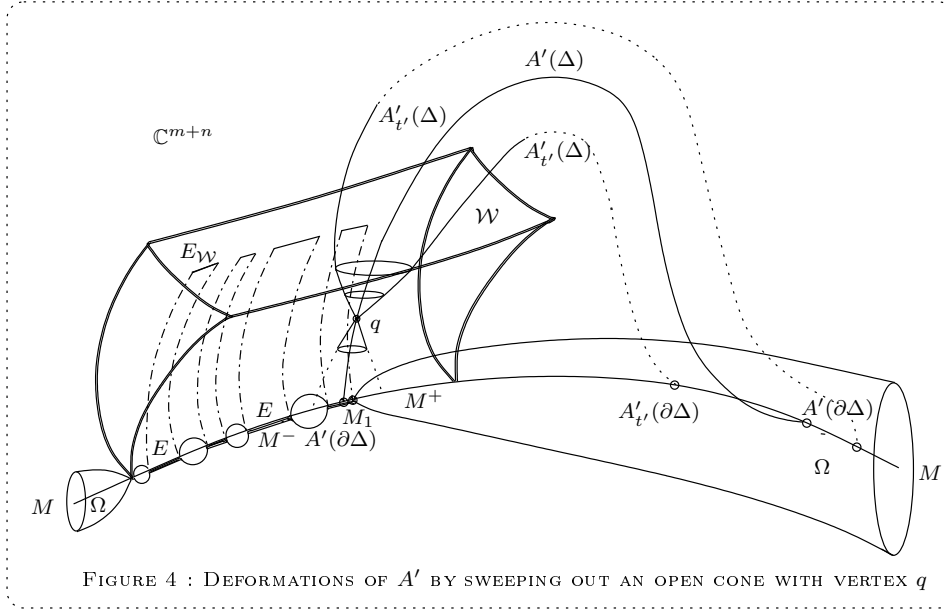
\put(585,-2931){\makebox(0,0)[lb]{\smash{\SetFigFont{8}{9.6}{\familydefault}{\mddefault}{\updefault}{\sc Figure 4: Deformations of $A'$ by sweeping out an open cone with vertex $q$}}}}
\put(1981,-2053){\makebox(0,0)[lb]{\smash{\SetFigFont{8}{9.6}{\familydefault}{\mddefault}{\updefault}$A'(\partial\Delta)$}}}
\put(2395,-1318){\makebox(0,0)[lb]{\smash{\SetFigFont{8}{9.6}{\familydefault}{\mddefault}{\updefault}$q$}}}
\put(5101,-1854){\makebox(0,0)[lb]{\smash{\SetFigFont{8}{9.6}{\familydefault}{\mddefault}{\updefault}$A'(\partial\Delta)$}}}
\put(3364,339){\makebox(0,0)[lb]{\smash{\SetFigFont{8}{9.6}{\familydefault}{\mddefault}{\updefault}$A'(\Delta)$}}}
\put(3055,-650){\makebox(0,0)[lb]{\smash{\SetFigFont{8}{9.6}{\familydefault}{\mddefault}{\updefault}$\mathcal{ W}$}}}
\put(1198,-870){\makebox(0,0)[lb]{\smash{\SetFigFont{8}{9.6}{\familydefault}{\mddefault}{\updefault}$E_{\mathcal{ W}}$}}}
\put(1019,-2110){\makebox(0,0)[lb]{\smash{\SetFigFont{8}{9.6}{\familydefault}{\mddefault}{\updefault}$E$}}}
\put(637,-2467){\makebox(0,0)[lb]{\smash{\SetFigFont{8}{9.6}{\familydefault}{\mddefault}{\updefault}$\Omega$}}}
\put(3953,-1873){\makebox(0,0)[lb]{\smash{\SetFigFont{8}{9.6}{\familydefault}{\mddefault}{\updefault}$A_{t'}'(\partial\Delta)$}}}
\put(3357,-238){\makebox(0,0)[lb]{\smash{\SetFigFont{8}{9.6}{\familydefault}{\mddefault}{\updefault}$A_{t'}'(\Delta)$}}}
\put(2446,-13){\makebox(0,0)[lb]{\smash{\SetFigFont{8}{9.6}{\familydefault}{\mddefault}{\updefault}$A_{t'}'(\Delta)$}}}
\end{picture}

\bigskip

In geometric terms, {\bf (4')} tells that $A'$ can be included in a
family $A_{t'}'$ of discs passing through $q$ which {\it sweeps out an
open cone with vertex in $q$}. Using some basic
differential geometric computations, the reader can easily check that
the geometric property {\bf (4')} implies {\bf (4)} after adding one
supplementary real parameter $t_{2m+2n-1}'$ corresponding to the
radius $\rho=\v\zeta\v$ of the disc.  Fortunately, for the needs of
Step~5 below, the essential foliation property stated in {\bf (4)}
will be valuable with an only $\mathcal{ C}^{1,\beta}$-smooth family
and, as stated in the end of {\bf (4)}, this family yields an
appropriate disc $A_{t'}'$ with empty intersection with the
singularity, namely $A_{t'}'(\{\rho_\eps e^{i\theta}:
\theta\in\R\})\cap E_{\mathcal{ W}}^{\rm nr} =\emptyset$.  Using this
Main Lemma~4.3, we can now accomplish the last step of the proof of
Theorem~3.1.

\smallskip
\noindent
{\bf Step 5: Removal of the point $q\in E_{\mathcal{ W}}^{\rm nr}$.}  
Let $A_{t'}'$ the family that we obtain by applying Main Lemma~\ref{tr}
to $A'$. According to the last sentence of Main Lemma~4.3, we may
choose $t'$ arbitrarily small and a positive radius $\rho_\eps>0$
sufficiently small so that the boundary of analytic subdisc
$A_{t'}'(\{\rho_\eps e^{i\theta}: \theta\in\R\})$ does not intersect
$E_{\mathcal{ W}}^{\rm nr}$. Let us denote such a disc $A_{t'}'$ simply
by $A''$ in the sequel.  Furthermore, we can assume that
$A''(\{\rho_\eps e^{i\theta}: \theta\in\R\})$ is contained in the
small ball $B_\eps:=\{|z-q|\leq \eps\}$ in which we shall localize an
application of the continuity principle ({\it see} {\sc Figure~5}).
Thus, it remains essentially to check that $F$ extends analytically to a
neighborhood of $q$ in $\C^{m+n}$ by constructing an analytic isotopy
of $A''$ in $(\mathcal{ W}\backslash E_{\mathcal{ W}}^{\rm nr})
\cup\Omega$ and by
applying the continuity principle.

One idea would be to translate a little bit in $\C^{m+n}$ the small
disc $A''(\{\rho e^{i\theta}: \rho\leq\rho_\eps, \,
\theta\in\R\})$. However, there is {\it a priori} no reason for which
such a small translated disc (which is of real dimension two) would
avoid the singularity $E_{\mathcal{ W}}^{\rm nr}$. Indeed, since we
only know that $H^{2m+2n-1}(E_{\mathcal{ W}}^{\rm nr})=0$, it is
impossible in general that a two-dimensional manifold avoids such a
``big'' set of Hausdorff codimension $1^{+0}$.

Of course, there is no surprise here: it is clear that functions which
are holomorphic in the domain $\mathcal{ W}\backslash E_{\mathcal{
W}}^{\rm nr}$ do {\it not} extend automatically through a set with
$H^{2m+2n-1}(E_{\mathcal{ W}}^{\rm nr})=0$, since for instance, such a
set $E_{\mathcal{ W}}^{\rm nr}$ might contain infinitely many complex
hypersurfaces, which are certainly not removable. So we really need to
consider the whole disc $A''$ and to include it into another family of
discs attached to $M\cup \Omega_1$ in order to produce an appropriate
analytic isotopy.

The good idea is to include $A''$ in a family $A_{\rho,s}''$ similar
to the one in Section~2 (with of course $A_{\rho_0,0}=A''$, but
without the unnecessary parameter $v$), since we already know that for
almost all $s\in U_s$, we can show as in Step~1 above that
$f$ (hence $F$ too) extends holomorphically to a neighborhood of
$A_{\rho_0,s}(\overline{\Delta})$ in $\C^{m+n}$.

To construct this family, we observe that $A''$ is not attached to
$M$, but as $A''$ can be chosen arbitrarily close in $\mathcal{
C}^{1,\beta}$-norm to the original disc $A$ attached to $M$, it
follows that $A''$ is certainly attached to some $\mathcal{
C}^{1,\beta}$-smooth manifold $M''$ close to $M$ which coincides with
$M$ except in a neighborhood of $A''(-1)$. Finally, the family
$A_{\rho,s}''$ is constructed as in Section~2 (but without the
parameter $v$, because in order to add the parameter
$v$ satisfying the second order
condition {\bf (4)} of Section~3, one would need $\mathcal{
C}^{2,\beta}$-smoothness of the disc).  By Tumanov's regularity
theorem [Tu2], this family is again of class $\mathcal{
C}^{1,\beta}$ for all $0<\beta<\alpha$.  Using properties {\bf (3)}
and {\bf (6)} and reasoning as in Step~1 of this Section~4 (continuity
principle), we deduce that the function $f$ of Theorem~3.1 extends
holomorphically to a neighborhood of $A_{\rho_0,s}(\overline{\Delta})$
in $\C^{m+n}$ for all $s\in U_s$, except those belonging to some
closed thin set $\mathcal{ S}$ with $H^{2m+n-2}(\mathcal{
S})=0$. Since $\mathcal{ S}$ does not locally disconnect $U_s$, such an
extension necessarily coincides with the extension $F$ in the
intersection of their domains.

In summary, by using the family $A_{\rho,s}''$, we have shown that 
for almost all $s$, {\it the function $F$ extends holomorphically to 
a neighborhood of $A_{\rho_0,s}''(\overline{\Delta})$}. We can therefore
apply the continuity principle to {\it remove the point $q$}.

Indeed, we remind that $A''=A_{\rho_0,0}''$ and that by construction the
small boundary $A_{\rho_0,0}''(\{\rho_\eps e^{i\theta}:
\theta\in\R\})$ which encircles $q$ does not intersect $E_{\mathcal{
W}}^{\rm nr}$. It is now clear that the usual continuity principle
along the family of small discs $A_{\rho_0,s}''(\{\rho e^{i\theta}:
\rho<\rho_\eps,\, \theta\in\R\})$ yields holomorphic extension of $F$
at $q$ ({\it see} again {\sc Figure~5}).

\bigskip
\hspace{-1cm}
\begin{picture}(0,0)%
\epsfig{file=figure5.pstex}%
\end{picture}%
\setlength{\unitlength}{3947sp}%
\begingroup\makeatletter\ifx\SetFigFont\undefined%
\gdef\SetFigFont#1#2#3#4#5{%
  \reset@font\fontsize{#1}{#2pt}%
  \fontfamily{#3}\fontseries{#4}\fontshape{#5}%
  \selectfont}%
\fi\endgroup%
\begin{picture}(6144,3054)(199,-2368)
\put(1939,-49){\makebox(0,0)[lb]{\smash{\SetFigFont{8}{9.6}{\familydefault}{\mddefault}{\updefault}$E_{\mathcal{ W}}$}}}
\put(4058,406){\makebox(0,0)[lb]{\smash{\SetFigFont{8}{9.6}{\familydefault}{\mddefault}{\updefault}$A''(\Delta)$}}}
\put(2746,-532){\makebox(0,0)[lb]{\smash{\SetFigFont{8}{9.6}{\familydefault}{\mddefault}{\updefault}$B_\varepsilon$}}}
\put(2626,-438){\makebox(0,0)[lb]{\smash{\SetFigFont{8}{9.6}{\familydefault}{\mddefault}{\updefault}$q$}}}
\put(5096,-1168){\makebox(0,0)[lb]{\smash{\SetFigFont{8}{9.6}{\familydefault}{\mddefault}{\updefault}$A''(\partial\Delta)$}}}
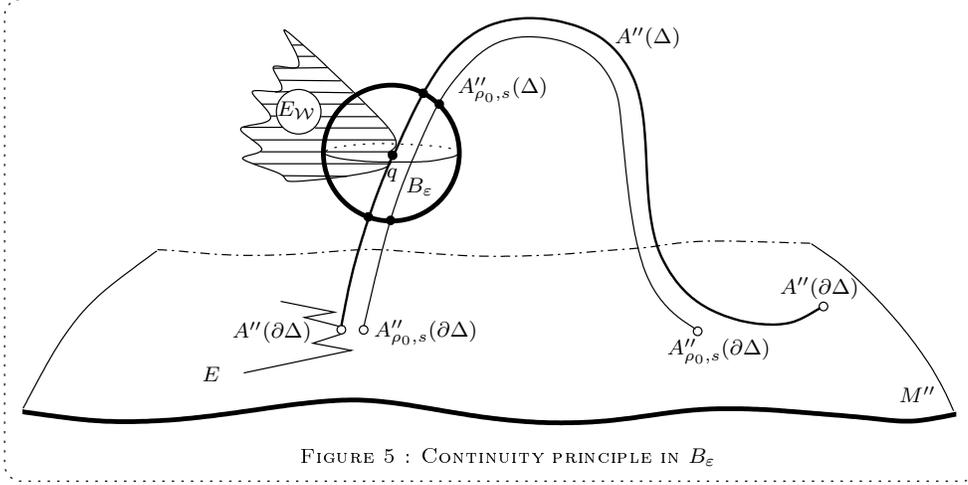
\put(2078,-2234){\makebox(0,0)[lb]{\smash{\SetFigFont{8}{9.6}{\familydefault}{\mddefault}{\updefault}{\sc Figure 5: Continuity principle in $B_\varepsilon$}}}}
\put(1467,-1721){\makebox(0,0)[lb]{\smash{\SetFigFont{8}{9.6}{\familydefault}{\mddefault}{\updefault}$E$}}}
\put(3071, 94){\makebox(0,0)[lb]{\smash{\SetFigFont{8}{9.6}{\familydefault}{\mddefault}{\updefault}$A_{\rho_0,s}''(\Delta)$}}}
\put(2547,-1436){\makebox(0,0)[lb]{\smash{\SetFigFont{8}{9.6}{\familydefault}{\mddefault}{\updefault}$A_{\rho_0,s}''(\partial\Delta)$}}}
\put(4387,-1552){\makebox(0,0)[lb]{\smash{\SetFigFont{8}{9.6}{\familydefault}{\mddefault}{\updefault}$A_{\rho_0,s}''(\partial\Delta)$}}}
\put(1657,-1446){\makebox(0,0)[lb]{\smash{\SetFigFont{8}{9.6}{\familydefault}{\mddefault}{\updefault}$A''(\partial\Delta)$}}}
\put(5843,-1852){\makebox(0,0)[lb]{\smash{\SetFigFont{8}{9.6}{\familydefault}{\mddefault}{\updefault}$M''$}}}
\end{picture}

\bigskip

Finally, the proof of Theorem~3.1 is complete modulo the proof of
Main Lemma~4.3, to which the remainder of the paper is devoted.

\section{Analytic discs attached to maximally real manifolds}

A crucial ingredient of the proof of Theorem~\ref{red} is the
description of a family of analytic discs which are close to the given
disc $A'$ of Main Lemma~4.3 and which are attached to a maximally 
real submanifold $R\subset M\cup \Omega_1$ (we shall construct
such an $R$ with $A'(\zeta)\in R$ for each $\zeta\in\partial\Delta$ in
Section~6 below).  This topic was developed by E.~Bedford--B.~Gaveau,
F.~Forstneri$\check{\rm c}$ in complex dimension 
two and generalized by J.~Globevnik
to higher dimensions. In this introductory section, we shall closely
follow [G1,2].

We need the solution of the following more general distribution
problem. Instead of a fixed maximally real submanifold $R$, we
consider a smooth family $R(\zeta)$, $\zeta\in\partial\Delta$, of
maximally real submanifolds of $\C^N$, $N\geq 2$, and we study the
discs attached to this family which are close to an attached disc $A'$
of reference, {\it i.e.}~fulfilling $A'(\zeta)\in R(\zeta)$, $\forall
\, \zeta\in\partial\Delta$. Let $\alpha>0$ be as in Theorem~1.1 and
let $0<\beta<\alpha$ be arbitrarily close to $\alpha$, as in Main
Lemma~4.3.

Concretely, the manifolds $R(\zeta)$ are given by defining
functions $r_j\in \mathcal{ C}^{2,\beta}(\partial\Delta\times B,
\R)$, $j=1,\dots,N$, where $B\subset\C^N$ is a small
open ball containing the origin, so that $r_j(\zeta,0)=0$ and
$\partial r_1(\zeta,p)\wedge\dots\wedge \partial r_N(\zeta,p)$
never vanishes for $\zeta\in\partial\Delta$ and $p\in B$. We would like
to mention that in [G1,~p.~289], the author considers the more general
regularity $r_j\in\mathcal{ C}^\beta(\partial\Delta,\mathcal{
C}^2(B))$, but that for us, the simpler smoothness category 
$\mathcal{ C}^{2,\beta}(\partial\Delta\times B,\R)$ will be enough.
Then we  represent
$$
R(\zeta):=\{p\in A(\zeta)+B: r_j(\zeta,p-A'(\zeta))=0,j=1,\dots,N\},
$$
which is a $\mathcal{C}^{2,\beta}$-smooth maximally real
manifold by the condition on
$\partial r_j$. We suppose the given reference disc $A'$ to be of
class $\mathcal{ C}^{2,\beta}$ up to the boundary. Following [G1], we
describe the family of nearby attached disc as a $\mathcal{
C}^{1,\beta}$-smooth submanifold of the space $\mathcal{
C}^{2,\beta}(\partial\Delta, \C^N)$, with a loss of smoothness.

\smallskip
\noindent 
{\it Remark.}
At first glance the transition from a fixed manifold to the family
$R(\zeta)$ may appear purely technical. Nevertheless it gives
in our application a decisive additional degree of freedom:
If we had to construct a fixed manifold $R$ containing the boundary
of our given disc $A'$, the boundary of $A'$ would prescribe 
one direction of $TR$. It will prove very convenient
to avoid this constraint by the transition to distributions
$R(\zeta)$ and this freedom will be used in an essential way in 
Section~6 below.

\smallskip
It turns out that the problem is governed by an $N$-tuple
$\kappa_1,\dots,\kappa_N\in\Z$ of coordinate independent partial
indices which are defined as follows.  As in [G1], we shall always
assume that the pull-back bundle $(A'|_{\partial\Delta})^*(TR(\zeta))$
is topologically trivial (this condition is dispensible, {\it
see}~[O]). For each $\zeta\in\partial\Delta$, let us denote by
$L(\zeta)$ the tangent space to $R(\zeta)$ at $A'(\zeta)$. Then there
is a $\mathcal{ C}^{1,\beta}$-smooth map $G:\partial\Delta\rightarrow
GL(N,\C)$ such that for each $\zeta\in\partial\Delta$ the columns of
$G$ are a (real) basis of $L(\zeta)$. By results of Plemelj and Vekua,
we can decompose the matrix function
$B(\zeta)=G(\zeta)\overline{G(\zeta)^{-1}},\zeta\in\partial\Delta$, as
$$
B(\zeta)=F^+(\zeta)\Lambda(\zeta)F^-(\zeta),
$$
with matrix functions 
\begin{eqnarray*}
F^+\in\mathcal{O}(\Delta,GL(N,\C))\cap 
\mathcal{ C}^{1,\beta}(\overline{\Delta},GL(N,\C)),& &\\
F^-\in\mathcal{O}(\C\backslash\overline{\Delta},GL(N,\C))\cap 
\mathcal{ C}^{1,\beta}(\C\backslash\Delta,GL(N,\C)),& \ &
\end{eqnarray*}
and where $\Lambda(\zeta)$ is the matrix with powers
$\zeta^{\kappa_j}$ on the diagonal and zero elsewhere. In [G1] it is
shown that the matrix $B(\zeta)$ depends only on the family 
of maximally real linear
space $L(\zeta)$ and that the $\kappa^j$ are unique up to
permutation. They are called the {\it partial indices} of $R$ along
$A'(\partial\Delta)$ and their sum $\kappa=\kappa_1+\dots+\kappa_N$
the {\it total index}. We stress that only $\kappa$ is a topological
invariant, in fact twice the winding number of $\det(G(\zeta))$ around
the origin. In the literature on symplectic topology, $\kappa$ is
called {\it Maslov index} of the loop $\zeta\mapsto L(\zeta)$.

Building on work of Forstneri$\check{\rm c}$ [F], Globevnik [G1,
Theorem~7.1] showed that the family of all analytic discs attached to
$R(\zeta)$ which are $\mathcal{ C}^{1,\beta}$-close to $A'$ is a
$\mathcal{ C}^{1,\beta}$-smooth submanifold of $\mathcal{O}(\Delta,\C^N)\cap
\mathcal{ C}^{2,\beta}(\overline{\Delta},\C^N)$ of dimension $\kappa+N$,
if all $\kappa_j$ are non-negative (by a result due to Oh [O], this
is even true if $\kappa_j\geq -1$ for all $j$). Furthermore the result
is stable with respect to small $\mathcal{ C}^{2,\beta}$-smooth 
deformations of $M$.

We shall need some specific ingredients of Globevnik's construction.
Since all our later arguments will exclude the appearance of odd
partial indices and since the expression of the square root matrix
$\sqrt{\Lambda}$ below is less complicated for even ones, we shall suppose
from now on that $\kappa_j=2m_j,j=1,\dots,N$.

Firstly one has to replace $G(\zeta)$ by another basis of
$L(\zeta)$ which extends holomorphically to
$\Delta$. By [G1, Lemma~5.1], there is a finer decomposition
$$
B(\zeta)=\Theta(\zeta)\Lambda(\zeta)\overline{\Theta(\zeta)^{-1}},
$$
where $\Theta\in\mathcal{O}(\Delta,GL(N,\C))\cap 
\mathcal{ C}^{1,\beta}(\overline{\Delta},GL(N,\C))$.
The substitute for $G(\zeta)$ is 
$$
\Theta(\zeta)\sqrt{\Lambda}(\zeta),
$$
where $\sqrt{\Lambda}(\zeta)$ denotes the matrix with $\zeta^{m_j}$ on the
diagonal. We denote by $X_j$ ($Y_j$) the columns of
$\Theta(\zeta)\sqrt{\Lambda}(\zeta)$ ($\sqrt{\Lambda}(\zeta)$)
respectively. One can verify that the $X_j(\zeta)$ span
$L(\zeta)$ ([G1, Theorem~5.1]).  Observe
$\Theta(\zeta)\sqrt{\Lambda}(\zeta)\in \mathcal{O}(\Delta,\C^N)\cap
\mathcal{ C}^{1,\beta}(\overline{\Delta},\C^N)$.

Secondly one studies variations of $A'|_{\partial\Delta}$ as a
function from $\partial\Delta$ to $\C^N$. Every nearby $\mathcal{
C}^{1,\beta}$-smooth (not necessarily 
holomorphic) variation is a disc close to $A'$ which can be written in
the form ([F, p.~20])
$$                                            
G(u,f)(\zeta)=\sum_{j=1}^N \, u_j(\zeta)\, X_j(\zeta)+i\, \sum_{j=1}^N \,
\{ f_j(\zeta)+i\, (T_0 f_j)(\zeta)\} \, X_j(\zeta),
$$
where $u_j,f_j\in\mathcal{C}^{1,\beta}(\partial\Delta,\R)$, are
uniquely determined by the variation.  Here $T_0$ denotes the harmonic
conjugation operator normalized at $\zeta=0$.  The condition
$G(u,f)(\zeta)\in R(\zeta)$, $\forall\, \zeta\in\partial \Delta$ is
equivalent to the validity of the system
$r_j(\zeta)(G(u,f)(\zeta))=0$, $1\leq j\leq N$.  The implicit function
theorem implies that this system can be solved for $f=\phi(u)$ for
$\mathcal{ C}^{1,\beta}$-small $u$ with a $\mathcal{
C}^{1,\beta}$-smooth mapping $\phi$ of Banach spaces $\mathcal{C}^{1,
\beta}(\partial\Delta,\R^N)\to \mathcal{C}^{1,
\beta}(\partial\Delta,\R^N)$. This follows from [G1, Theorem~6.1] by
an application of the implicit function theorem in Banach spaces,
except concerning the $\mathcal{ C}^{1,\beta}$-smoothness, which, in
our situation, is more direct and elementary than in [G1], since we
have supposed that $r_j\in\mathcal{ C}^{2,\beta}(\partial\Delta\times
B,\R)$.

Finally one has to determine for which choices of $u$ the function
$G(u,\phi(u))$ extends holomorphically to $\Delta$. Writing
$$
G(u,f)(\zeta)=\Theta(\zeta)\, \sum_{j=1}^N \left\{u_j(\zeta)
 +i\, [f_j(\zeta)+i\, (T_0 f_j)(\zeta)]\right\} \, Y_j,
$$
we see that $G(u,\phi(u))$ extends holomorphically, if
and only if
\def\theequation{5.1}
\begin{equation}
\Theta^{-1}(\zeta)\, G(u,\phi(u))(\zeta)= \sum_{j=1}^N \{u_j(\zeta)
 +i\, [\phi(u)_j(\zeta)+i\, (T_0 \,\phi(u)_j)(\zeta)]\}\, Y_j
\end{equation}
extends, {\it i.e.}~if and only if the function
$\zeta\mapsto\sum_{j=1}^N u_j(\zeta)\,Y_j(\zeta)$
extends. One can compute ([G1, p.~301]) that this is precisely the
case, if $h_j(\zeta)=Y^{-1}(\zeta)\,u_j(\zeta)$ has polynomial components of
the form
\def\theequation{5.2}
\begin{equation}
\aligned
h_j(\zeta) = & \ t_1^j+i\,t_2^j + (t_3^j+i\,t_4^j)\, \zeta + \dots + 
  (t_{\kappa_{j-1}}^j+i\,t_{\kappa_j}^j)\, \zeta^{m_j-1} 
  + t_{\kappa_{j+1}}^j\,\zeta^{m_j} \\
    & + (t_{\kappa_{j-1}}^j-i\,t_{\kappa_j}^j)\,\zeta^{m_j+1}  
  + \dots + (t_3^j-i\,t_4^j)\,\zeta^{\kappa_j-1}
+ (t_1^j-i\,t_2^j)\,\zeta^{\kappa_j},
\endaligned
\end{equation}
where all $t_k^j$ are real. In total we get
$\kappa_j+1$ real parameters for the choice of $h_j$ and hence
$\kappa+N$ parameters for our local family of discs attached to
$R(\zeta)$.

\section{Proof of Theorem~3.1, part II}

In this section we provide the final part of the proof of
Theorem~\ref{red}, namely Main Lemma~4.3, which relies essentially on
global properties of analytic discs. The disc $A'$ of
Main Lemma~4.3 need not be attached to $M$ but since it is
close to $A$ in $\mathcal{ C}^{2,\beta}$-norm, it is certainly
attached to a nearby manifold $M'$ of class $\mathcal{ C}^{2,\beta}$
which coincides with $M$ except in $\Omega_1$. The idea is now to
first embed $A'(\partial\Delta)$ into a maximal real submanifold of
$M\cup \Omega_1$ whose partial indices are easy to determine. Then we
shall explain how to increase the partial indices separately by
twisting $R$ around $A'(\partial\Delta)$ inside $\Omega_1$. The
families of attached discs get richer with increasing indices and will
eventually contain the required discs $A_{t'}'$ as a subfamily. We
divide the proof in four essential steps.

\smallskip
\noindent
{\bf Step 1: Construction of a first maximally real manifold
$R_1$.} Let $h'$ be a defining function of $M'$ as in~(3.1). Then $A'$
is the solution of a Bishop equation
$$
Y'=T_1(h'(W',Y'))+y_0,
$$
where $W'\in \mathcal{ C}^{2,\beta}$ is the $w$-component of $A'$ and
$y_0\in\R^n$ is close to $0$. Recall that by construction, $W'(\zeta)$
is close to the $w$-component $(\rho_0-\rho_0\zeta,0,\dots, 0)$ of
the disc $A$ defined in~(3.3). Let $A_{u_*,y}'$ be the discs defined by
the perturbed equation
$$
Y_{u_*,y}'=T_1(h'(W_{u_*,y}'+(0,u_*),Y_{u_*,y}'))+y_0+y,
$$
where $u_*:=(u_2,\dots,u_m)$ is close to $0$ and $y\in\R^n$ is close to
$0$. We have $A_{0,0}'=A'$. Since $A$ defined by~(3.3) and hence
also $A'$ are by construction almost parallel to the $w_1$-axis, the
union
$$
R_1:=\bigcup_{u_*,y}\,A_{u_*,y}'(\partial\Delta)
$$
is a maximally real manifold of class $\mathcal{ C}^{2,\beta}$
contained in $M'$ and containing $A'(\partial\Delta)$. 
The explicit construction of $R_1$ allows
an easy determination of the partial indices.

\begin{lem}\label{ind}
The partial indices of $R_1$ with respect to $A'(\partial\Delta)$
are $2,0,\dots,0$.
\end{lem}

\proof
We begin by constructing $N=m+n$ holomorphic vector fields along
$A'(\partial\Delta)$ which generate (over $\R$) the tangent bundle of
$R_1$. We denote $\zeta=e^{i\theta}\in\partial \Delta$ and define
first $G_1(\zeta):=[\partial A'(e^{i\theta})/\partial \theta]$ as the
push-forward of $\partial/\partial\theta$. Next, we put
$$
\left\{
\aligned
G_k(\zeta):= & \ 
[\partial A_{u_*,0}'(\zeta)/\partial u_k]|_{u_*=0}, \ \ \
\text{\rm  for} \ k=2,\dots,m,\\
G_k(\zeta):= & \
[\partial A_{0,y}'(\zeta)/\partial
y_{k-m}]|_{y=0}, \ \ \ \text{\rm  for} \
k=m+1,\dots,N.
\endaligned\right.
$$  
For $k=2,\dots,N$, $G_k$ is the
uniform limit of pointwise holomorphic difference quotients and
therefore holomorphic itself. As $A_{u_*,y}'$ depends $\mathcal{
C}^{2,\beta}$-smoothly on parameters, we obtain $G_k\in \mathcal{
C}^{1,\beta}(\overline{\Delta},\C^N)$, $k=2,\dots,N$.

By [G1, Proposition 10.2], the maximal number of
linearly independent holomorphically extendable sections equals the
number of non-negative partial indices. Hence we deduce that all
$\kappa_j$ are non-negative.

Furthermore it is easy to see that the total index $\kappa$, which is
twice the winding number of ${\rm det}\, G|_{\partial\Delta}$ around
$0$, equals $2$.  Indeed, $A'$ is almost parallel to the $w_1$ axis,
the direction in which $G_1$ has winding number $1$, and the vector
fields $G_2,\dots,G_N$ have a topologically trivial behaviour in the
remaining directions. This heuristic argument can be made precise in
the following way. One easily can smoothly deform the complex
coordinates $z_j,w_k$ to (non-holomorphic coordinates) in which the
matrix $G(\zeta)$ gets diagonal with diagonal entries
$\zeta,1,\dots,1$.  In the deformed coordinates the winding number of
the determinant is obviously $1$, and this remains unchanged when
deforming back to the standard coordinates.

In summary the only possible constellations for the partial indices
are $2,0,\dots$, $0$ and $1,1,0,\dots,0$. But [G1, Proposition 10.1]
excludes the second case as $\partial A'/\partial\theta$ does not
vanish on $\overline{\Delta}$, which completes the proof.
\endproof

\smallskip
\noindent
{\bf Step 2: Gluing $R_1$ with a family of maximally real planes.}
Our goal is to twist the manifold $R_1$ many times around the boundary
of $A'$ in the small neighborhood $\Omega_1$ of $A'(-1)$ in order to
{\it increase} its partial indices. Since it is rather easy to
increase partial indices when a disc is attached to a family of {\it
linear} maximally real subspaces of $\C^N$ (using Lemma~6.3 below, {\it
see} the reasonings just after the proof), we aim to glue $R_1$ with
its family of tangent planes $T_{A'(\zeta)}R_1$ for $\zeta$ near
$-1$. Before proceeding, we have to take care of a regularity question:
the family $\zeta\mapsto T_{A'(\zeta)}R_1$ being only of class
$\mathcal{ C}^{1,\beta}$, some preliminary regularizations are
necessary. We remind that by Lemma~4.2 {\bf (b)}, the disc $A'$ is
real analytic near $\zeta=-1$. This choice of smoothness is very
adapted to our purpose. Indeed, using cut-off functions and the
Weierstrass approximation theorem, we can construct a $\mathcal{
C}^{2,\beta}$-smooth maximally real manifold $R_2$ to which $A'$ is
still attached and which is also {\it real analytic} in a neighborhood
of $\{A'(\zeta): \v\zeta+1\v<\eps/2\}$. Of course, this can 
be done with $\v\v R_2-R_1\v\v_{\mathcal{ C}^{2,\beta}}$ being arbitrarily
small, so the partial indices of $A'$ with respect to $R_2$ 
are still equal to $(2,0,\dots,0)$. 

Using real analyticity, we can now glue $R_2$ with its family of
maximally real tangent planes $T_{A'(\zeta)}R_2$ for
$\v\zeta+1\v<\eps/4$ in a {\it smooth way} as follows. After localization 
near $A'(-1)$ using a cut-off function, the gluing problem is reduced to the
following statement.

\begin{lem}
Let $R$ be small real analytic maximally real submanifold of $\C^N$,
let $p\in R$ and let $\gamma(s)$, $s\in (-\eps,\eps)$, be a real
analytic curve in $R$ passing through $p$. Then there exist smooth
functions $r_j(s,z)\in\mathcal{ C}^\infty((-\eps,\eps)\times B,\R)$,
for $j=1,\dots,N$, where $B$ is a small open ball centered at the
origin in $\C^N$, such that
\begin{itemize}
\item[{\bf (1)}]
$r_j(s,\gamma(s))\equiv 0$.
\item[{\bf (2)}]
$r_j(s,z)\equiv r_j(z)\equiv$ the defining functions of $R$ for 
$\v s\v\geq \eps/2$.
\item[{\bf (3)}]
For all $s$ with $\v s\v\leq\eps/4$, the set $\{z\in\C^N:
r_j(s,z)=0, j=1,\dots,N\}$ coincides with the tangent
space of $R$ at $\gamma(s)$.
\end{itemize}
\end{lem}

\proof
Choosing coordinates $(z_1,\dots,z_N)$ vanishing at $p$, we can
asssume that $R$ is given by $r_j(z):=y_j-\varphi_j(x)=0$ with 
$\varphi_j(0)=0$ and $d\varphi_j(0)=0$, and that
$\gamma_j(s)=x_j(s)+i y_j(s)$, where $y_j(s):=\varphi_j(x(s))$. Let
$\chi(s)$ be a $\mathcal{ C}^\infty$-smooth cut-off function
satisfying $\chi(s)\equiv 0$ for $\v s\v\leq\eps/4$ and $\chi(s)\equiv
1$ for $\v s\v\geq\eps/2$. We choose for $r_j(s,z)$ the following functions:
$$
y_j-y_j(s)-\sum_{k=1}^N
{\partial\varphi_j\over\partial x_k}(x(s))\,[x_k-x_k(s)]-
\chi(s)\left[\sum_{K\in\N^N,\v K\v\geq 2}
{\partial_x^K\varphi_j(x(s))\over K\,!}\,[x-x(s)]^K\right].
$$
Clearly, the $r_j$ are $\mathcal{ C}^\infty$-smooth
and {\bf (3)} holds.
As $\varphi_j$ is real analytic in a neighborhood
of $\gamma$, property {\bf (2)} holds by Taylor's formula.
\endproof

In summary, we have shown that we can attach $A'$ to some $\mathcal{
C}^{2,\beta}$-smooth family $(R_3(\zeta))_{\zeta\in\partial\Delta}$ of
maximally real submanifolds such that $R_3(\zeta)$ coincides with
$R_2$ for $\v \zeta+1\v\geq\eps/2$ and such that $R_3(\zeta)$
coincides with the maximally real plane $T_{A'(\zeta)}R_2$, for
$\v\zeta+1\v\leq\eps/4$. Clearly, the partial indices of $A'$ with
respect to the family $R_3(\zeta)$ are still equal to
$(2,0,\dots,0)$.

\smallskip
\noindent
{\bf Step 3: Increasing partial indices.}
This step is the crucial  one in our argumentation.
Recall that the partial indices are defined in terms of vector fields
along $A'(\partial\Delta)$. 
In the previous section we have described how to select
distinguished vector fields $X_k$ as the columns of
$\Theta\sqrt{\Lambda}$, where $\Theta,\Lambda$ were associated to a
decomposition of $G_3(\zeta)\,\overline{G_3(\zeta)}^{-1}$, 
where the columns of 
the matrix $G_3(\zeta)$ span $T_{A'(\zeta)}\,R_3(\zeta)$. 
Our method is to modify the vector fields
$X_k$ by replacing them by products $g_k X_k$ with the boundary values
of certain holomorphic functions $g_k$. 
It turns out that the indices can be read from
properties of the $g_k$. Here is how the $g_k$ are constructed.

For convenience in the following lemma, we shall represent
$\partial\Delta$ by the real closed interval $[-\pi,\pi]$ where $\pi$ is
identified with $-\pi$.

\begin{lem}
For every small $\eps>0$, every integer $\ell\in\N$, there exists a
holomorphic function $h\in\mathcal{ O}(\Delta)\cap\mathcal{
C}^\infty(\overline{\Delta})$ such that
\begin{itemize}
\item[{\bf (1)}] 
$h(\zeta)\neq 0$ for all $\zeta\in\overline{\Delta}$.
\item[{\bf (2)}] 
The function $g(\zeta):=\zeta^\ell\,h(\zeta)$ is real-valued over
$\{e^{i\theta}: \v\theta\v\leq \pi-\eps/8\}$.
\end{itemize}
It follows that the winding number 
of $g\v_{\partial\Delta}$ around $0\in\Delta$ is equal to $\ell$.
\end{lem}

\proof
Let $v(\zeta)$ be an arbitrary $\mathcal{ C}^\infty$-smooth
$2\pi$-periodic extension to $\R$ of the linear function $-\ell\theta$
defined on $[-\pi+\eps/8,\pi-\eps/8]$.  Let $T_0$ be the harmonic
conjugate operator satisfying $(T_0\,u)(0)=0$ for every $u\in
L^2(\partial\Delta)$.  Since $T_0$ is a bounded operator of the
$\mathcal{ C}^{k,\alpha}$ spaces of norm equal to $1$, the function
$T_0\,v$ is $\mathcal{ C}^\infty$-smooth over $\partial\Delta$.  It
suffices to set $h:=\exp(-T_0\,v+iv)$. Indeed,
$$
\zeta^\ell\,h(\zeta)=e^{i\ell\theta}\, e^{-T_0\,v+iv}
$$
is real for $\theta\in [-\pi+\eps/8,\pi-\eps/8]$, as desired.
\endproof

Let $L(\zeta)$ denote the tangent space $T_{A'(\zeta)}\,R_3(\zeta)$ and let 
$X_k(\zeta)$ be $\mathcal{ C}^{1,\beta}$-smooth vector fields
as the columns of the matrix $\Theta\sqrt\Lambda$ constructed
in Section~5 above. We remind that the $X_k(\zeta)$
span $L(\zeta)$. Further, as the partial indices of
$A'(\zeta)$ with respect to $R_3(\zeta)$ are $(2,0,\dots,0)$, 
we have 
$$
\Theta(\zeta)\sqrt{\Lambda}(\zeta)=
(\zeta\,\Theta_1(\zeta),\Theta_2(\zeta),\dots,\Theta_N(\zeta)).
$$ 
Now, let us choose an arbitrary collection of nonnegative integers
$\ell_1,\ell_2,\dots,\ell_N$ and associated functions
$g_{\ell_1}(\zeta)=\zeta^{\ell_1}\,h_{\ell_1}(\zeta)$, $\dots$,
$g_{\ell_N}(\zeta)=\zeta^{\ell_N}\,h_{\ell_N}(\zeta)$ satisfying 
{\bf (1)} and
{\bf (2)} of Lemma~6.3 above. With these functions, we define a new family of
maximally real manifolds to which $A'(\zeta)$ is still attached as
follows:

\begin{itemize}
\item[{\bf (a)}]
For $\v\theta\v\leq \pi-\eps/8$, $R_4(\zeta)\equiv R_3(\zeta)$.
\item[{\bf (b)}]
For $\v\theta\v\geq\pi-\eps/8$, namely for $\zeta$ close to $-1$,
$$
R_4(\zeta):=\text{\rm span}_\R\,(
\zeta\,g_{\ell_1}(\zeta)\,\Theta_1(\zeta),
g_{\ell_2}(\zeta)\,\Theta_2(\zeta),\dots,
g_{\ell_N}(\zeta)\,\Theta_N(\zeta)).
$$
\end{itemize}
It is important to notice that this definition yields a true
$\mathcal{ C}^{2,\beta}$-smooth family of maximally real manifolds, 
thanks to the fact that the family $R_3(\zeta)$ is already a family of
{\it real linear spaces} for $\v \zeta+1\v \leq \eps/4$, by 
construction. Interestingly, the partial indices have increased:

\begin{lem}
The partial indices of $R_4(\zeta)$ along $\partial A$ are equal to 
$2+2\ell_1,2\ell_2,\dots$, $2\ell_N$.
\end{lem}

\proof
By construction, since the functions $g_{\ell_j}(\zeta)$ are real-valued
over $\{e^{i\theta}: \v\theta\v\leq \pi-\eps/8\}$, the tangent space
$T_{A'(\zeta)}R_4(\zeta)$ is spanned {\it for all}
$\zeta\in\partial\Delta$ by the $N$ vectors
$$
\zeta\,g_{\ell_1}(\zeta)\,\Theta_1(\zeta),\,
g_{\ell_2}(\zeta)\,\Theta_2(\zeta),\,\dots,\,
g_{\ell_N}(\zeta)\,\Theta_N(\zeta),
$$
which form together a $N\times N$ matrix which we will denote by
$G_4(\zeta)$.  By Section~5, we can read directly from the matrix
identity
$$
\aligned
G_4(\zeta)\,& \overline{G_4(\zeta)}^{-1} =
(h_{\ell_1}(\zeta)\,\Theta_1(\zeta),\dots,
h_{\ell_N}(\zeta)\,\Theta_N(\zeta))\,\times\\
&
\times\text{\rm diag}\,(\zeta^{2+2\ell_1},\zeta^{2\ell_2},\dots,
\zeta^{2\ell_N})\,\times \,
\overline{(h_{\ell_1}(\zeta)\,\Theta_1(\zeta),\dots,
h_{\ell_N}(\zeta)\,\Theta_N(\zeta))}^{-1}
\endaligned
$$
that the partial indices of $A'(\zeta)$ with respect to $R_4(\zeta)$
are equal to $(2+2\ell_1,2\ell_2,\dots$, $2\ell_N)$, as stated.
\endproof

\smallskip
\noindent
{\bf Step 4: Construction of the family $A_{t'}'$.} 
Now, as we need not very large partial indices, we choose $\ell_1=1$,
$\ell_2=2,\dots,\ell_N=2$, so the partial indices are simply
$(4,4,\dots,4)$.
Moreover, the matrix $Y(\zeta)$ is equal to the
diagonal matrix ${\rm diag}\,(\zeta^2,\zeta^2,\dots,\zeta^2)$. 
Concerning the $5N$ parameters
$(t_1^j,t_2^j,t_3^j,t_4^j,t_5^j)$ appearing in equation~(5.2) above,
we even choose $t_1^j=t_2^j=t_5^j=0$. Then by the result of Globevnik,
we thus obtain a family of discs depending on the $2N$-dimensional
real parameter $t:=(t_3^j+i\,t_4^j)_{1\leq j\leq N}$.
The functions $h_j$ and $u_j$ defined in Section~5 are thus equal to
$$
\left\{
\aligned
h_j(t,\zeta):= & \ (t_3^j+i\,t_4^j)\,\zeta+(t_3^j-i\,t_4^j)\,\zeta^3,\\
u_j(t,\zeta):= & \ (t_3^j+i\,t_4^j)\,\bar\zeta+(t_3^j-i\,t_4^j)\,\zeta.
\endaligned\right.
$$
It remains to explain how we can extract the desired family $A_{t'}'$
by reducing this $(2m+2n)$-dimensional parameter space to some
of dimension $(2m+2n-1)$ such that property {\bf (4)}
of Main Lemma~4.3 is satisfied. Let us denote by $h_t$ and
$u_t$ the maps $\zeta\mapsto h(t,\zeta)$ and
$\zeta\mapsto u(t,\zeta)$. By equation~(5.1), we have
$$
G(u_t,\phi(u_t))(\zeta)=
\Theta(\zeta)\,\sum_{j=1}^N\left\{
u_j(t,\zeta)+i\,\left[\phi(u_t)_j(\zeta)+i\,T_0\,\phi(u_t)_j(\zeta)\right]
\right\}\,Y_j(\zeta).
$$
and by Section~5, the $\mathcal{ C}^{1,\beta}$-smooth discs 
$$
A_t'(\zeta):=A'(\zeta)+G(u_t,\phi(u_t))(\zeta)
$$
are attached to $R_4(\zeta)$. By [G1, p.~299 top], the differential
of $\phi$ at $0$ is null: $D_u\phi(0)=0$. It follows that
$$
{\partial\over\partial t}\,\left[
\Theta(\zeta)\,\sum_{j=1}^N\,
i\,\left[
\phi(u_t)_j(\zeta)+i\,T_0\,\phi(u_t)_j(\zeta)
\right]\,Y_j(\zeta)
\right]_{t=0}\,\equiv 0.
$$
So on the one hand, we can compute for $j=1,\dots,N$
\def\theequation{6.1}
\begin{equation}
\left\{
\aligned
\left[
{\partial A_t'\over\partial t_3^j}\right]_{t=0}= & \
\rho e^{i\theta}\,
\Theta(0)\,(0,\dots,0,1,0,\dots,0)+{\rm O}(\rho^2),\\
\left[
{\partial A_t'\over\partial t_4^j}\right]_{t=0}= & \
\rho e^{i\theta}\,
\Theta(0)\,(0,\dots,0,i,0,\dots,0)+{\rm O}(\rho^2),\\
\endaligned\right.
\end{equation}
where ${\rm O}(\rho^2)$ denotes a holomorphic disc in $\mathcal{
O}(\Delta,\C^N) \cap\mathcal{ C}^{0,\beta}(\overline{\Delta},\C^N)$
vanishing up to order one at $0$. For $\rho>0$ small enough and $\theta$
arbitrary, it follows that these
$2m+2n$ vectors span $\C^{m+n}$.  On the other hand, we compute
\def\theequation{6.2}
\begin{equation}
\left[
{\partial A_t'\over\partial \theta}\right]_{t=0}=
\rho e^{i\theta}\,
\Theta(0)\,(ia_1,\dots,ia_N)+{\rm O}(\rho^2),
\end{equation}
where the constants $a_j$ are defined by $A'(\zeta)=(a_1\zeta,\dots,
a_N\zeta)+{\rm O}(\zeta^2)$ and do not all vanish (since $A'$ is
an embedding). 

Let us choose a $(2m+2n-1)$-dimensional real plane $H$ which is
supplementary to $\R\,\Theta(0)\,(ia_1,\dots,ia_N)$ in $\C^{m+n}$.
Using~(6.1), we can choose a $(2m+2n-1)$-dimensional
real linear subspace $T'\subset\R^{2m+2n}$ and $\rho_\eps$ small enough 
such that, after restricting the family $A_t'$ with $t'\in T'$,
the $(2m+2n-1)$ vectors $[\partial A_{t'}'/\partial t_j']_{t'=0}$, 
$j=1,\dots,2m+2n-1$, 
are linearly independent with the vector~(6.2) for all $\zeta\in\Delta$
of the form $\zeta=\rho_\eps e^{i\theta}$. It follows that the mapping
$$
(e^{i\theta},t')\mapsto A_{t'}'(\rho_\eps e^{i\theta})
$$
is a local embedding of the circle $\partial\Delta$ times a small
neighborhood of the origin in $\R^{2m+2n-1}$, from which we see that
the foliation property {\bf (4)} of Main Lemma~4.3 holds.

This completes the proof of Step 4, the proof of Main Lemma~4.3,
the proof of Theorem~3.1 and the proof of Theorem~1.1.

\end{document}